\documentclass[12pt]{article}
\usepackage{ae,german,amsthm,amssymb,a4,amsmath,makeidx,epsfig,amsfonts,amscd}
\usepackage{psfrag}
\usepackage{makeidx}

\makeindex

\newtheorem{stz}{Satz}
\newtheorem{lmm}{Lemma}
\newtheorem{krl}{Korollar}
\theoremstyle{definition}
\newtheorem{bmr}{Bemerkung}
\theoremstyle{definition}
\newtheorem{dfn}{Definition}

\DeclareMathOperator{\supp}{supp}
\DeclareMathOperator{\de}{d}
\DeclareMathOperator{\dist}{dist}
\DeclareMathOperator{\diam}{diam}
\DeclareMathOperator{\con}{con}
\DeclareMathOperator{\re}{Re}
\DeclareMathOperator{\im}{Im}
\DeclareMathOperator{\e}{e}
\DeclareMathOperator{\GL}{GL}

\begin{document}

\title{\bf Kohomologie von Garben nichtabelscher Gruppen}
\author{Diplomarbeit}
\date{Humboldt-Universität zu Berlin \\ Mathematisch-Naturwissenschaftliche Fakultät II \\ Institut für Mathematik}
\maketitle

\noindent \\[5cm] eingereicht von: Katrin Kaden \\ $ $ \\ geb.: am 27.07.76 in: Halle/Saale \\ $ $ \\ Betreuer: Prof. Dr. J. Leiterer \\ $ $ \\ Berlin, den 9. September 2004

\leavevmode
\newpage

\noindent \\[8cm] 
An dieser Stelle möchte ich Herrn Professor Dr. J. Leiterer für das interessante Thema, seine intensive Betreuung und seine vielen wertvollen Hinweise danken. Weiterhin danke ich meiner Familie und meinen Freunden für ihre Unterstützung und für ihr Verständnis, was mit zum Gelingen dieser Diplomarbeit beigetragen hat.

\newpage 

\tableofcontents

\newpage

\section{Einführung}
\label{heinrich}

\noindent Kohomologie-Methoden werden in der komplexen Analysis seit etwa 1950 intensiv studiert, zum Beispiel von J. Leray, H. Cartan, J. P. Serre und K. Stein (siehe \cite{gr2}, Seite XVII). Die dabei untersuchten mathematischen Strukturen sind auf topologischen Räumen durch gewisse Mengen stetig lokal definierter Abbildungen, auch Abbildungskeime genannt, erklärt. Mathematisch exakt lassen sich diese Strukturen als Garben definieren, was auf J. Leray zurückgeht (siehe \cite{lr1} und \cite{lr2}). In Abhängigkeit von Eigenschaften der betrachteten Abbildungen beschäftigen wir uns besonders mit Garben von stetigen und holomorphen Abbildungskeimen, wobei die Bilder dieser Abbildungen in der multiplikativen topologischen Gruppe einer Banachalgebra liegen sollen, die nicht notwendig abelsch ist. Dadurch verallgemeinern wir den in \cite{gr2} und \cite{gun} ausführlich behandelten Fall von matrixwertigen Abbildungen, betrachten im holomorphen Fall allerdings nur eine komplexe Veränderliche. Es werden deshalb Garben über nicht notwendigerweise abelschen Gruppen eingeführt. Leider ist die Kohomologietheorie im nichtabelschen Fall nur ansatzweise ausgearbeitet. Das liegt unter anderem daran, dass die Gruppenstruktur der Kohomologiemenge $H^1$ bereits verloren geht. 

\noindent Im Abschnitt 2.1 behandeln wir die grundlegende Theorie der Garben über nicht notwendig abelschen Gruppen, und führen die erste Kohomologiemenge $H^1$ ein. Es ergeben sich Schwierigkeiten bei der Definition höherer Kohomologiemengen $H^q$, $q=2,3$, worauf in dieser Arbeit nicht eingegangen wird. 

\noindent Die grunglegende Absicht der Arbeit besteht in der Beantwortung der Frage, unter welchen Umständen die Kohomologiemenge $H^1$ trivial ist. Diesbezügliche Zerfällungssätze formulieren und beweisen wir sowohl im stetigen Fall (Abschnitt 3.1) unter relativ allgemeinen Forderungen an $X$, als auch im holomorphen Fall (Abschnitt 3.2) unter spezielleren Voraussetzungen. In den Beweisen werden Faktorisierungssätze verwendet. Diese sind im holomorphen Fall analog den Cousinschen beziehungsweise Cartanschen Heftungslemmata (siehe zum Beispiel \cite{gr2}, Kapitel III) und werden für unsere Fälle speziell formuliert und bewiesen. Die Faktorisierungssätze benutzen im stetigen Fall einen Fortsetzungssatz, der mit Hilfe des Fortsetzungssatzes von Dugundji (siehe \cite{lem} oder \cite{gro}) bewiesen wird, und im holomorphen Fall einen Approximationssatz und den Satz von Graves (siehe \cite{lan}). Beim Übergang vom paarweisen Schnitt zu endlichen Überdeckungen war es zweckmä\3ig, sogenannte $\cal F$-Ketten für den stetigen Fall und $\cal F$-Felder besonders für den holomorphen Fall zu definieren (Abschnitt 2.2), wobei $\cal F$ eine beliebige Garbe über einer nicht notwendig abelschen Gruppe ist. Mit Satz \ref{harry-lucia} enthält der Abschnitt 2.1 eine relativ allgemeine Aussage, deren Wert erst im Beweis des letzten Satzes der Arbeit beim Nachweis des Zerfällungssatzes im holomorphen Fall deutlich wird. Bei der Anwendung dieses Satzes ergibt sich schlie\3lich  durch einen Grenzübergang, dass die betrachtete Kohomologiemenge $H^1$ trivial ist. Man kann die in den Hauptsätzen der Arbeit (Satz \ref{horst-petra}, Satz \ref{hector-gabi}) verwendeten Aussagen schematisch folgenderma\3en darstellen:

\begin{figure}
 \begin{center}
  \includegraphics*[width=15cm]{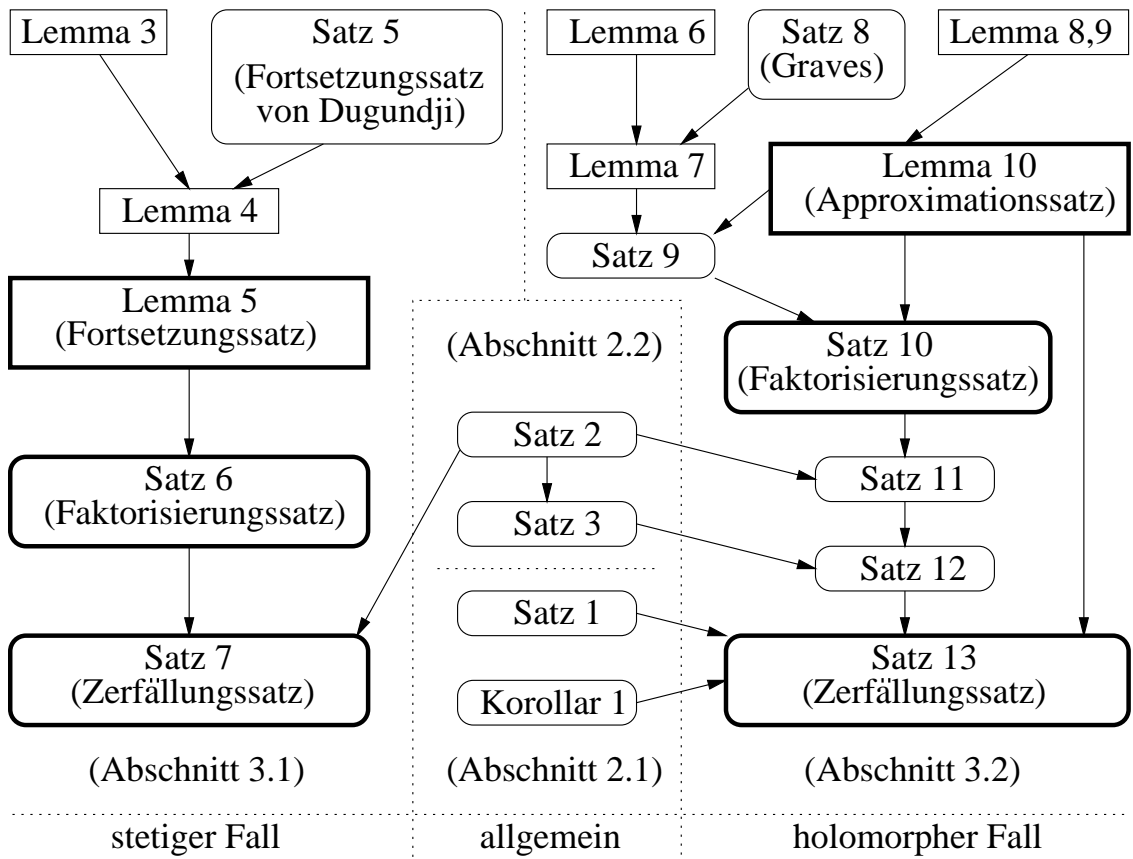}
 \end{center} 
 \caption{\label{fig:s} Logische Abhängigkeiten der Sätze und Hilfssätze}
\end{figure}

\newpage

\section{Grundlegende Definitionen und Sätze}
\label{harry}

\subsection{Kohomologiemengen}
\label{har}

\noindent Für die nachfolgenden Definitionen vergleiche man \cite{hir}, \cite{kul} und \cite{gin}.
\begin{dfn} Eine {\sl Garbe}\index{Garbe} $\cal G$ von nicht notwendigerweise abelschen Gruppen über $X$ (kurz: eine Garbe über X) ist ein Tripel $(S,\pi,X)$ mit folgenden Eigenschaften:
\begin{enumerate}
\item $S$ und $X$ sind topologische Räume und $\pi:S\rightarrow X$ ist eine surjektive, stetige Abbildung, welche {\sl Projektion}\index{Projektion einer Garbe} der Garbe genannt wird.
\item Jeder Punkt $s\in S$ besitzt eine offene Umgebung $N\subset S$, so dass $\pi|_N$ ein Homöomorphismus von $N$ auf        $\pi(N)$ ist. 
\item Die Menge ${\cal G}_x:=\pi^{-1}(x)$, der sogenannte {\sl Halm}\index{Halm}, ist für jedes $x\in X$ eine Gruppe.
\item Die durch $(s,t)\mapsto st^{-1}$ für jedes $x\in X$ definierte Abbildung von ${\cal G}_x\times {\cal G}_x$ in ${\cal G}_x$ ist stetig.
\end{enumerate}
\end{dfn}

\psfrag{X}{$X$}
\psfrag{S}{$S$}
\psfrag{s}{$s$}
\psfrag{Gx}{${\mathcal G}_x$}
\psfrag{x}{$x$}
\psfrag{pi}{$\pi$}
\psfrag{N}{$N$}
\psfrag{piN}{$\pi(N)$}
\begin{figure}[htb]
 \begin{center}
  \includegraphics*[width=8cm]{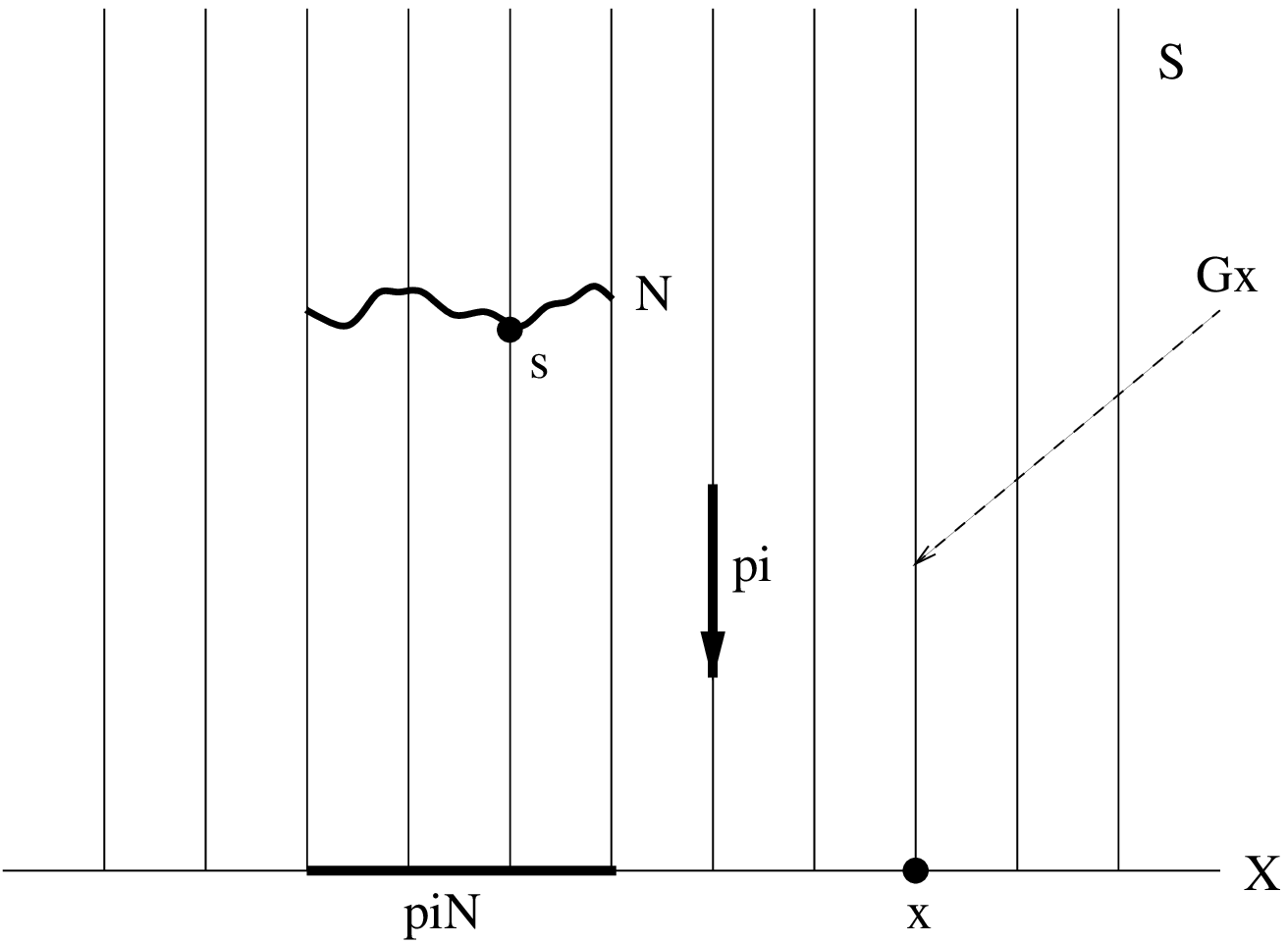}
 \end{center} 
 \caption{\label{fig:DefGarbe} Definition einer Garbe}
\end{figure}

\begin{dfn} Ein {\sl Schnitt}\index{Schnitt} in einer Garbe $\cal G$ über einer offenen Menge $U\subset X$ ist eine stetige Abbildung $f:U\rightarrow S$, für die $\pi\circ f=\mathrm{id}$ gilt. Die Menge aller Schnitte in $\cal G$ über $U$ ist eine Gruppe (mit der Operation $(fg)(x):=f(x)g(x), x\in U$, $f,g$ Schnitte über $U$) und wird mit ${\cal G}(U)$ bezeichnet. Das Einselement der Gruppe ist der {\sl Einsschnitt}\index{Einsschnitt} $x\mapsto 1_x\in{\cal G}_x$, $x\in U$ und wird mit $1$ bezeichnet. Inverse Elemente $f^{-1}$ werden in ${\cal G}(U)$ durch $(f^{-1})(x):=f(x)^{-1}, x\in U$ definiert.\footnote{$f^{-1}$ bezeichnet also nicht die Umkehrabbildung von $f$.}
\end{dfn}

\noindent Ein einfaches Beispiel einer Garbe ${\cal G}=(S,\pi,X)$ von nichtabelschen Gruppen ist folgenderma\3en definiert: Sei $S:=X\times\GL(r),\ r\geq 2$ für einen topologischen Raum $X$, wobei $\GL(r)\subset{\mathbb C}^{r\times r}$ die Gruppe der invertierbaren $r\times r$-Matrizen ist. Dabei betrachte man auf $\GL(r)$ die diskrete Topologie, womit insbesondere einpunktige Teilmengen von $\GL(r)$ offen sind. Auf $S$ betrachte man die Produkttopologie. Es sei $\pi:(x,A)\mapsto x$ für $x\in X$ und $A\in \GL(r)$ die Projektion des direkten Produktes. Für ein beliebiges Element $s=(x,A)\in S$ sei $V$ eine offene Umgebung von $x$ in $X$, so dass $N:=V\times\{A\}\subset S$ eine offene Umgebung von $s$ bezüglich der Produkttopologie ist. Dann ist $\pi|_N$ offenbar ein Homöomorphismus von $N$ auf $V$. In diesem Beispiel einer Garbe $\cal G$ kann jeder Halm ${\cal G}_x, x\in X$ mit $\GL(r)$ identifiziert werden. Über einer offenen Menge $U\subset X$ erfüllen die Schnitte $f:U\rightarrow U\times\GL(r)$, definiert durch $f(x):=(x,A_x),\ A_x\in\GL(r)$, die Bedingung $\pi\circ f=\mathrm{id}$, und sie sind stetig. \\[0,5cm]

\noindent Kohomologie{\sl gruppen} $H^q(X,{\cal G})$ lassen sich für $q\geq1$ im nichtabelschen Fall nicht definieren, jedoch kann für $q=1$ eine Kohomologie{\sl menge} $H^1(X,{\cal G})$ mit ausgezeichnetem Element definiert werden. 

\noindent Sei dazu ${\cal U}=(U_i)_{i\in I}$ eine offene Überdeckung von $X$. Zur Vereinfachung benutzen wir die Bezeichnungen $U_{ij}:=U_i\cap U_j$ und $U_{ijk}:=U_i\cap U_j\cap U_k$ für alle $i,j,k\in I$. Eine Familie $h=(h_i)_{i\in I}$ mit $h_i\in {\cal G}(U_i)$ für alle $i\in I$ hei\3t {\sl $0$-Kokette}\index{0-Kokette} über ${\cal U}$. Es sei ${\cal C}^0({\cal U,G})$ die Menge aller $0$-Koketten über $\cal U$. Analog nennen wir eine Familie $f=(f_{ij})_{i,j\in I}$ mit $f_{ij}\in {\cal G}(U_{ij})$ für alle $i,j\in I$ $1$-{\sl Kokette}\index{1-Kokette} und bezeichnen mit ${\cal C}^1({\cal U,G})$ die Menge aller $1$-Koketten über $\cal U$. Die Mengen ${\cal C}^0({\cal U,G})$ und ${\cal C}^1({\cal U,G})$  sind Gruppen, deren Operation durch das Produkt von Schnitten erklärt wird. Eine $0$-Kokette $h=(h_i)_{i\in I}$ hei\3t $0$-{\sl Kozyklus}\index{0-Kozyklus}, wenn $$h_i^{-1}h_j=1$$ auf $U_{ij}$ für alle $i,j\in I$ gilt. Es sei ${\cal Z}^0({\cal U},{\cal G})$ die Menge aller $0$-Kozyklen über $\cal U$. Eine $1$-Kokette $f=(f_{ij})_{i,j\in I}$ mit $$f_{ij}f_{jk}=f_{ik}$$ auf $U_{ijk}$ für alle $i,j,k\in I$ hei\3t $1$-{\sl Kozyklus}\index{1-Kozyklus} über $\cal U$. Aus dieser Definition folgt sofort, dass dann $f_{ii}$ der Einsschnitt über $U_i$ ist und $f_{ji}=f_{ij}^{-1}$ gilt. Die Menge aller $1$-Kozyklen über $\cal U$ bezeichnet man mit ${\cal Z}^1(\cal U,\cal G)$. 

\noindent Man kann ${\cal Z}^0({\cal U},{\cal G})$ mit der Gruppe ${\cal G}(X)=:H^0(X,{\cal G})$ identifizieren, denn wegen $(h_i)_{i\in I}\in{\cal Z}^0({\cal U},{\cal G})$ mit $h_i^{-1}h_j=1$, also $h_i=h_j$ auf $U_{ij}$ für alle $i,j\in I$, ist genau ein Schnitt $h\in{\cal G}(X)$ durch $h|_{U_i}=h_i,\ i\in I$ definiert.

\noindent Die Menge ${\cal Z}^1({\cal U},{\cal G})$ ist im Allgemeinen keine Untergruppe der Kokettengruppe ${\cal C}^1({\cal U},{\cal G})$: Für den Nachweis sei beispielsweise ${\cal U}:=(U_1,U_2)$ mit $U_1\cap U_2\neq\emptyset$ und $f,g\in{\cal Z}^1({\cal U},{\cal G})$, wobei $f_{12}$ und $g_{12}$ nicht vertauschbar seien. Wegen $(f\cdot g)_{21}=f_{21}\cdot g_{21}=(f_{12})^{-1}(g_{12})^{-1}\neq(f_{12}\cdot g_{12})^{-1}=((f\cdot g)_{12})^{-1}$ ist $f\cdot g\notin{\cal Z}^1({\cal U},{\cal G})$.

\noindent Für $f=(f_{ij})_{i,j\in I}\in{\cal Z}^1({\cal U},{\cal G})$ und $h=(h_i)_{i\in I}\in{\cal C}^0({\cal U,G})$ bezeichnet man mit $h\Box f$ den $1$-Kozyklus, der durch $$(h\Box f)_{ij}:=h_i^{-1}f_{ij}h_j$$ auf $U_{ij}$ für alle $i,j\in I$ definiert ist.

\begin{dfn} Zwei $1$-Kozyklen $f,g\in {\cal Z}^1(\cal U,\cal G)$ hei\3en {\sl äquivalent}\index{äquivalente 1-Kozyklen}, in Zeichen $f\sim g$, wenn ein $h\in {\cal C}^0({\cal U,G})$ existiert, so dass $$g=h\Box f$$  gilt. Die zu dieser Äquivalenzrelation gehörige Menge von Äquivalenzklassen bezeichnet man mit $H^1(\cal U,\cal G)$. Für $f\in{\cal Z}^1({\cal U},{\cal G})$ benutzt man die Bezeichnung $[f]$ für die zu $f$ gehörige Klasse aus $H^1(\cal U,\cal G)$.
\end{dfn}

\noindent Für jede offene Überdeckung ${\cal U}=(U_i)_{i\in I}$ von $X$ bezeichnet man mit $1\in H^1(\cal U,\cal G)$ die Äquivalenzklasse aller $f=(f_{ij})_{i,j\in I}\in {\cal Z}^1(\cal U,\cal G)$, für die eine $0$-Kokette $h=(h_i)_{i\in I}\in {\cal C}^0({\cal U,G})$ existiert, so dass $$f_{ij}=h_i^{-1}h_j$$ auf $U_{ij}$ für alle $i,j\in I$ gilt, also $f=h\Box 1$. Man sagt in diesem Falle auch, dass die 1-Kozyklen $f$ {\sl zerfallen}\index{zerfallende 1-Kozyklen}. Das Element $1$ hei\3t {\sl neutrales Element}\index{neutrales Element} oder {\sl Einselement}\index{Einselement} von $H^1(\cal U,\cal G)$.

\noindent Eine offene Überdeckung ${\cal V}=(V_j)_{j\in J}$ von $X$ hei\3t {\sl Verfeinerung}\index{Verfeinerung, allgemein} von ${\cal U}=(U_i)_{i\in I}$, wenn es eine Abbildung $\varphi:J\rightarrow I$ gibt mit $V_j\subset U_{\varphi j}$ für alle $j\in J$. Folglich induziert $\varphi$ eine Abbildung $$\varphi^*:{\cal Z}^1({\cal U},{\cal G})\rightarrow {\cal Z}^1({\cal V},{\cal G})\:,\ \text{die durch}\ \:(\varphi^*f)_{ij}=f_{\varphi i,\varphi j}|_{V_{ij}}\:,i,j\in J\:,\ f\in{\cal Z}^1({\cal U},{\cal G})$$ definiert ist. Diese induziert dann weiter eine Abbildung $$\varphi_{\cal V}^{\cal U}:H^1({\cal U},{\cal G})\rightarrow H^1(\cal V,\cal G)\:.$$ Im folgenden Lemma \ref{harry-susi} wird gezeigt, dass die Abbildung $\varphi_{\cal V}^{\cal U}$ nicht von der Wahl von $\varphi:J\rightarrow I$ abhängt.

\begin{lmm}\label{harry-susi} Sei $\psi:J\rightarrow I$ eine weitere Abbildung mit $V_j\subset U_{\psi j}$ für alle $j\in J$ und sei $f\in {\cal Z}^1(\cal U,\cal G)$. Dann gilt $[\varphi^*f]=[\psi^*f]$.
\end{lmm}

\begin{proof}[{\bf Beweis}] Man definiert eine $0$-Kokette $(h_j)_{j\in J}\in {\cal C}^0({\cal V,G})$ durch $$h_j:=f_{\varphi j,\psi j}|_{V_j}\in {\cal G}(V_j)$$ für alle $j\in J$. Nun gilt: $$(\psi^*f)_{ij}=f_{\psi i,\psi j}=f_{\psi i,\varphi i}f_{\varphi i,\varphi j}f_{\varphi j,\psi j}=h_i^{-1}(\varphi^*f)_{ij}h_j$$ auf $V_{ij}$ für alle $i,j\in J$, woraus aufgrund der Definition von $H^1(\cal V,\cal G)$ die Behauptung folgt. 
\end{proof}

\noindent Weiterhin zeigt man: 

\begin{lmm}\label{harry-peggy} Die Abbildung $\varphi_{\cal V}^{\cal U}:H^1({\cal U},{\cal G})\rightarrow H^1(\cal V,\cal G)$ ist injektiv.
\end{lmm}

\begin{proof}[{\bf Beweis}] Seien $\Phi_1,\Phi_2\in H^1(\cal U,\cal G)$ und $f,g\in {\cal Z}^1(\cal U,\cal G)$ mit $\Phi_1=[f]$ und $\Phi_2=[g]$. Es gelte $[\varphi^*f]=[\varphi^*g]$ 
und man zeige $[f]=[g]$. Nach der Voraussetzung gibt es eine $0$-Kokette $(h_j)_{j\in J}\in {\cal C}^0({\cal V,G})$ mit $$(\varphi^*g)_{ij}=h_i^{-1}(\varphi^*f)_{ij}h_j\:,\ \text{also}\ \:g_{\varphi i,\varphi j}=h_i^{-1}f_{\varphi i,\varphi j}h_j$$ auf $V_{ij}$ für alle $i,j\in J$. Weiterhin gilt: $$g_{\varphi i,k}g_{k,\varphi j}=h_i^{-1}f_{\varphi i,k}f_{k,\varphi j}h_j\:,\ \text{d.h.}\ \:f_{k,\varphi i}h_ig_{\varphi i,k}=f_{k,\varphi j}h_jg_{\varphi j,k}$$ auf $U_k\cap V_{ij}$ für alle $i,j\in J$ und $k\in I$. Folglich ist $l_k\in {\cal G}(U_k)$ für $k\in I$ durch $$l_k|_{U_k\cap V_j}:=f_{k,\varphi j}h_jg_{\varphi j,k}|_{U_k\cap V_j}$$ für alle $j\in J$ mit $U_k\cap V_j\neq\emptyset$ wohldefiniert, also $(l_k)_{k\in I}\in {\cal C}^0({\cal U,G})$ und es gilt: $$l_m^{-1}f_{mn}l_n=g_{m,\varphi j}h_j^{-1}f_{\varphi j,m}f_{mn}f_{n,\varphi j}h_jg_{\varphi j,n}=g_{m,\varphi j}h_j^{-1}f_{\varphi j,\varphi j}h_jg_{\varphi j,n}=g_{mn}$$ auf $U_{mn}\cap V_j$ für alle $m,n\in I$ und $j\in J$, woraus $\Phi_1=[f]=[g]=\Phi_2$ folgt und die Injektivität gezeigt ist.
\end{proof}

\noindent Als Folgerung daraus ergibt sich unmittelbar:

\begin{krl}\label{harry-clothilde} Aus $H^1({\cal V},{\cal G})=1$ folgt $H^1({\cal U},{\cal G})=1$, wenn $\cal V$ eine Verfeinerung von $\cal U$ ist.
\end{krl}

\noindent Es sei $\tilde{H}^1(X,\cal G)$ die disjunkte Vereinigung aller $H^1(\cal U,\cal G)$, wobei $\cal U$ alle offenen Überdeckungen von $X$ durchläuft. (Dabei soll zugelassen werden, dass eine offene Menge mehrmals in einer Überdeckung vorkommt. Um logische Schwierigkeiten zu vermeiden, betrachte man nur Indexmengen, die Teilmengen einer hinreichend gro\3 gewählten Menge sind.) 

\noindent Zwei Elemente $\Phi_1\in H^1(\cal U,\cal G)$ und $\Phi_2\in H^1(\cal V,\cal G)$ (wobei $\cal U$ und $\cal V$ irgendwelche offenen Überdeckungen von $X$ sind) hei\3en {\sl äquivalent}\index{äquivalente Klassen von 1-Kozyklen} $\Phi_1\sim\Phi_2$, wenn es eine offene Überdeckung $\cal W$ von $X$ gibt, die sowohl $\cal U$ als auch $\cal V$ verfeinert, so dass $$\varphi_{\cal W}^{\cal U}\Phi_1=\varphi_{\cal W}^{\cal V}\Phi_2$$ gilt. Die {\sl Kohomologiemenge}\index{Kohomologiemenge} $H^1(X,\cal G)$ ist dann definiert als die Menge aller Äquivalenzklassen von $\tilde{H}^1(X,\cal G)$ bezüglich dieser Äquivalenzrelation. 

\noindent Dabei sind die neutralen Elemente von $H^1({\cal U},{\cal G})$, wobei $\cal U$ alle offenen Überdeckungen von $X$ durchläuft, offenbar paarweise äquivalent. Die dazugehörige Klasse in $H^1(X,{\cal G})$ wird ebenfalls mit $1$ bezeichnet. Besteht die Kohomologiemenge nur aus dem Einselement, so sagt man, dass sie {\sl trivial}\index{triviale Kohomologiemenge} ist.

\noindent Sei $\cal U$ eine offene Überdeckung von $X$. Für $\Phi\in H^1({\cal U},{\cal G})$ bezeichnen wir die durch $\Phi$ in $H^1(X,{\cal G})$ definierte Klasse mit $[\Phi]$. 

\noindent Es ist komplizierter, höhere Kohomologiemengen, zum Beispiel $H^2(X,\cal G)$ und $H^3(X,\cal G)$, zu definieren, weshalb wir diese hier nicht betrachten (vergleiche dazu \cite{ser}, Seite 46 und siehe beispielsweise \cite{gir}).

\noindent Seien $Y\subset X$ eine offene Teilmenge und ${\cal U}=(U_i)_{i\in I}$ eine offene Überdeckung von $X$. Für $f\in{\cal Z}^1({\cal U},{\cal G})$ definiert man durch $$f|_Y:=(f_{ij}|_{U_{ij}\cap Y})_{i,j\in I}$$ ein Element $f|_Y\in{\cal Z}^1({\cal U}|_Y,{\cal G})$, wobei ${\cal U}|_Y:=(U_i\cap Y)_{i\in I}$ ist und nennt $f|_Y$ {\sl Einschränkung}\index{Einschränkung von 1-Kozyklen} des 1-Kozyklus $f$ auf $Y$. Für $\Phi\in H^1({\cal U},{\cal G})$ wählt man $f\in{\cal Z}^1({\cal U},{\cal G})$ mit $[f]=\Phi$ und setzt $$\Phi|_Y:=[f|_Y]\in H^1({\cal U}|_Y,{\cal G})\:.$$ Man sieht leicht, dass diese Definition nicht von der Wahl von $f$ abhängt. Ist schlie\3lich $\Psi\in H^1(X,{\cal G})$, so wählt man eine offene Überdeckung $\cal U$ von $X$ sowie ein $\Phi\in H^1({\cal U},{\cal G})$ mit $[\Phi]=\Psi$ und setzt $$\Psi|_Y:=[\Phi|_Y]\in H^1(Y,{\cal G})\:.$$ Man sieht wieder leicht, dass diese Definition nicht von der Wahl von $\cal U$ und $\Phi$ abhängt. Wir nennen wieder $\Phi|_Y$ und $\Psi|_Y$ {\sl Einschränkungen}\index{Einschränkung von Klassen von 1-Kozyklen} von $\Phi$ beziehungsweise $\Psi$ auf $Y$.

\begin{stz}\label{harry-lucia} Es seien $\Psi\in H^1(X,\cal G)$ und ${\cal U}=(U_i)_{i\in I}$ eine offene Überdeckung von $X$. Gilt $\Psi|_{U_n}=1$ für alle $n\in I$, so existiert ein $\Phi\in H^1({\cal U},{\cal G})$ mit $[\Phi]=\Psi$.
\end{stz}

\begin{proof}[{\bf Beweis}] Laut Definition von $H^1(X,\cal G)$ gibt es eine offene Überdeckung ${\cal V}=(V_j)_{j\in J}$ von $X$, so dass ein $f=(f_{ij})_{i,j\in J}\in{\cal Z}^1(\cal V,\cal G)$ existiert mit $[[f]]=\Psi$. Ohne Beschränkung der Allgemeinheit sei $\cal V$ eine Verfeinerung von $\cal U$, d.h. es gibt eine Abbildung $\varphi:J\rightarrow I$ mit $V_j\subset U_{\varphi j}$ für alle $j\in J$.

\noindent Wegen $\Psi|_{U_n}=1\ \forall n\in I$ existieren nun  $0$-Koketten $(h_{jn})_{j\in J}\in {\cal C}^0({\cal V}|_{U_n},\cal G)$ mit $$f_{ij}=h_{in}^{-1}h_{jn}$$ auf $V_{ij}\cap U_n$ für alle $i,j\in J$ und für alle $n\in I$. Wir wählen $m,n\in I$ beliebig und definieren $$g_{mn}:=h_{jm}h_{jn}^{-1}$$ auf $U_{mn}\cap V_j$ mit $j\in J$. Dabei sind die so definierten $g_{mn}$ von der Wahl der $j\in J$ unabhängig, denn es gilt:  $$h_{in}^{-1}h_{jn}=f_{ij}=h_{im}^{-1}h_{jm}\:,\ \text{also}\ \:h_{im}h_{in}^{-1}=h_{jm}h_{jn}^{-1}$$
auf $U_{mn}\cap V_{ij}$, $i\in J$. Weiterhin gilt: $$g_{ml}g_{ln}=h_{jm}h_{jl}^{-1}h_{jl}h_{jn}^{-1}=g_{mn}$$ auf $U_{mln}\cap V_j$ für alle $m,l,n \in I$ und $j\in J$, so dass $g=(g_{mn})_{m,n\in I}$ ein $1$-Kozyklus über $\cal U$ ist, also $[g]\in H^1(\cal U,\cal G)$. Schlie\3lich gilt $[\varphi^*g]=[f]\in H^1(\cal V,\cal G)$ wegen $$(\varphi^*g)_{ij}=g_{\varphi i,\varphi j}=h_{j,\varphi i}h_{j,\varphi j}^{-1}=h_{i,\varphi i}h_{i,\varphi i}^{-1}h_{j,\varphi i}\ h_{j,\varphi j}^{-1}=h_{i,\varphi i}f_{ij}|_{U_{\varphi i}}h_{j,\varphi j}^{-1}=h_i^{-1}f_{ij}h_j$$ auf $V_{ij}\subset U_{{\varphi i},{\varphi j}}$ für alle $i,j\in J$ mit der $0$-Kokette $(h_k)_{k\in J}\in {\cal C}^0({\cal V},\cal G)$, definiert durch $$h_k:=h_{k,\varphi k}^{-1}\in {\cal G}(V_k)$$ für alle $k\in J$. Damit ist $[g]\sim[f]$, also $[g]\in\Psi$, und folglich liefert das Element $\Phi:=[g]\in H^1(\cal U,\cal G)$ das Gewünschte. 
\end{proof}

\noindent Mit Hilfe von Satz \ref{harry-lucia} und Lemma \ref{harry-peggy} ergibt sich leicht:

\begin{krl} Es sei ${\cal U}=(U_i)_{i\in I}$ eine offene Überdeckung von $X$, so dass $H^1(U_n,{\cal G})$ für alle $n\in I$ trivial ist. Dann kann man $H^1({\cal U},{\cal G})$ mit $H^1(X,{\cal G})$ identifizieren.
\end{krl}

\subsection{${\cal F}$-Ketten und ${\cal F}$-Felder}
\label{ry}

\noindent Für die nächsten Sätze wird folgende Begriffsbildung benötigt.

\begin{dfn} Seien $\cal F$ eine Garbe über dem topologischen Raum $X$ und weiterhin $U_1,\ldots,U_n$, $n\in\mathbb N$ offene Teilmengen von $X$.\footnote{Dabei ist $\mathbb N$ die Menge der natürlichen Zahlen ohne Null.}
\begin{enumerate}
\item $(U_1,U_2)$ hei\3t $\cal F$-{\sl Paar}\index{$\cal F$-Paar}, falls für jedes $f\in {\cal F}(U_1\cap U_2)$ Elemente $f_1\in {\cal F}(U_1)$ und $f_2\in {\cal F}(U_2)$ existieren mit $f=f^{-1}_1f_2$ auf $U_1\cap U_2$.
\item Wir nennen ${\cal U}=(U_i)_{1\leq i\leq n}$ eine $\cal F$-{\sl Kette}\index{$\cal F$-Kette}, falls $(U_1\cup\ldots\cup U_i,U_{i+1})$ für jedes $i\in\{1,\ldots,n-1\}$ ein $\cal F$-Paar ist.
\end{enumerate}
\end{dfn}

\begin{stz}\label{harry-kathi} Für eine $\cal F$-Kette ${\cal U}$ gilt $H^1({\cal U,F})=1$.
\end{stz}

\begin{proof}[{\bf Beweis}] Der Beweis erfolgt durch vollständige Induktion. Sei ${\cal U}=(U_i)_{1\leq i\leq n}$ mit $n\in\mathbb N$ gegeben. Für $n=1$ ist die Aussage trivial. Nun seien $n>1$, $V:=U_1\cup\ldots\cup U_{n-1}$ und ${\cal V}:=(U_i)_{1\leq i\leq n-1}$. Es sei $f=(f_{ij})_{1\leq i,j\leq n}\in {\cal Z}^1({\cal U,F})$ beliebig gewählt. Dann gilt $\hat{f}:=(f_{ij})_{1\leq i,j\leq n-1}\in {\cal Z}^1({\cal V,F})$, und wegen $H^1({\cal V,F})=1$ gibt es $h_i\in {\cal F}(U_i)$, $i=1,\ldots,n-1$ mit $$f_{ij}=h^{-1}_ih_j$$ auf $U_i\cap U_j$ für alle $1\leq i,j\leq n-1$. Es sei $W:=V\cap U_n$. Jedes $x\in W$ liegt in einem $U_i$ $(i<n)$ und wir setzen $$g(x):=h_i(x)f_{in}(x)\:.$$ Dadurch ist $g$ auf $W$ wohldefiniert, denn auf $U_i\cap U_j$ $(i,j<n)$ gilt: $$h_if_{in}=h_if_{ij}f_{jn}=h_ih^{-1}_ih_jf_{jn}=h_jf_{jn}\:.$$ Da $\cal U$ eine $\cal F$-Kette ist, ist $(V,U_n)$ ein $\cal F$-Paar und für $g\in {\cal F}(W)$ gibt es ein $h\in {\cal F}(V)$ und ein $g_n\in {\cal F}(U_n)$ mit $$g=h^{-1}g_n$$ auf $W$. Für $1\leq i\leq n-1$ setzen wir $$g_i:=hh_i\in {\cal F}(U_i)\:.$$ 
\begin{eqnarray*}
\text{Sind}\ 1\leq i,j\leq n-1,\ \ \text{so gilt} & g^{-1}_ig_j=h^{-1}_ih^{-1}hh_j=h^{-1}_ih_j=f_{ij} & \text{auf}\ \ U_i\cap U_j.\\ 
\text{Ist}\ 1\leq\ i\ \ \leq n-1,\ \ \text{so gilt} & g^{-1}_ig_n=h^{-1}_ih^{-1}g_n\ =h^{-1}_ig\ =f_{in} & \text{auf}\ \ U_i\cap U_n. 
\end{eqnarray*}
Damit gilt $[f]=1$, und da $f$ beliebig gewählt war, ist also $H^1({\cal U,F})=1$.
\end{proof}

\noindent Sind ${\cal U}=(U_i)_{i\in I}$ und ${\cal V}=(V_j)_{j\in J}$ zwei Familien
von Teilmengen von $X$ mit disjunkten Indexmengen $I$ und $J$, dann sei
$${\cal U}\cup{\cal V}:=(W_k)_{k\in I\cup J}$$
die {\sl Vereinigung}\index{Vereinigung von Mengenfamilien} von ${\cal U}$ und ${\cal V}$, wobei $W_k:=U_k$ für
$k\in I$ und $W_k:=V_k$ für $k\in J$ gilt.

\begin{dfn} Seien ${\cal U}^j=(U_i^j)_{1\leq i\leq n_j},\ n_j\in{\mathbb N}$ für alle $j\in\{1,\ldots,m\},\ m\in\mathbb N$ ${\cal F}$-Ketten. Unter einem ${\cal F}$-{\sl Feld}\index{$\cal F$-Feld} verstehen wir die Vereinigung $${\cal U}:={\cal U}^1\cup\ldots\cup{\cal U}^m\:,$$ wobei ${\cal V}:=(U^j)_{1,\leq j\leq m}$ für $U^j:=\bigcup\limits_{1\leq i\leq n_j}U_i^j$ ebenfalls eine ${\cal F}$-Kette ist.
\end{dfn}

\begin{stz}\label{harry-conny} Für ein ${\cal F}$-Feld ${\cal U}$ gilt $H^1({\cal U},{\cal F})=1$.
\end{stz}

\begin{proof}[{\bf Beweis}] Es sei ${\cal U}:={\cal U}^1\cup\ldots\cup{\cal U}^m,\ m\in\mathbb N$ ein ${\cal F}$-Feld. Der Nachweis erfolgt mit vollständiger Induktion nach $m$. Für $m=1$ besteht ${\cal U}$ aus einer ${\cal F}$-Kette, und die Behauptung folgt direkt mit Satz \ref{harry-kathi}.

\noindent Nun seien $m>1$, ${\cal V}:={\cal U}^1\cup\ldots\cup{\cal U}^{m-1}=:(V_k)_{k\in K}$ und ${\cal U}^m=:(U_l)_{l\in I_m}$ für disjunkte Indexmengen $K$ und $I_m$ . Als Induktionsvoraussetzung gelte $H^1({\cal V},{\cal F})=1$. Es soll gezeigt werden, dass $H^1({\cal U},{\cal F})=1$ gilt. 

\noindent Es seien $V:=\bigcup\limits_{k\in K}V_k\ \text{und}\ U:=\bigcup\limits_{l\in I_m}U_l\ $. Dann ist $(V,U)$ nach Voraussetzung ein ${\cal F}$-Paar, und wir setzen $W:=V\cap U$.

\noindent Nun sei $f\in{\cal Z}^1({\cal U},{\cal F})$ beliebig gewählt, und es ist zu zeigen, dass ein $g\in{\cal C}^0({\cal U},{\cal F})$ existiert mit $f=g\Box 1$. Der 1-Kozyklus $f$ setzt sich folgenderma\3en zusammen: 
\begin{enumerate}
  \item $f_{ik}$ mit $i,k\in K$ als Schnitte in $V$,
	\item $f_{jl}$ mit $j,l\in I_m$ als Schnitte in $U$ und 
  \item $f_{kl}$ mit $k\in K,l\in I_m$ als Schnitte in $W$. 
\end{enumerate}  
Es sei $f'=(f_{ik})_{i,k\in K}\in{\cal Z}^1({\cal V},{\cal F})$. Nach der Induktionsvoraussetzung
gibt es ein $h'=(h_k)_{k\in K}\in{\cal C}^0({\cal V},{\cal F})$ mit $f'=h'\Box 1$.

\noindent Weiter sei $f''=(f_{jl})_{j,l\in I_m}\in{\cal Z}^1({\cal U}^m,{\cal F})$. Da ${\cal U}^m$ eine ${\cal F}$-Kette ist, gibt es nach Satz \ref{harry-kathi} ein $h''=(h_l)_{l\in I_m}\in{\cal C}^0({\cal U}^m,{\cal F})$ mit $f''=h''\Box 1$. Wir setzen $$h:=h_kf_{kl}h_l^{-1}$$ auf $V_k\cap U_l$ für $k\in K,l\in I_m$, wobei $h$ auf $W=\bigcup\limits_{k\in K,l\in I_m}(V_k\cap U_l)$ wegen $$h_kf_{kl}h_l^{-1}=h_kf_{ki}f_{ij}f_{jl}h_l^{-1}=h_kh_k^{-1}h_if_{ij}h_j^{-1}h_lh_l^{-1}=h_if_{ij}h_j^{-1}$$ auf $(V_i\cap V_k)\cap(U_j\cap U_l)$ für alle $i,k\in K, j,l\in I_m$ wohldefiniert ist, womit $h\in{\cal F}(W)$ gilt.
 
\noindent Da $(V,U)$ ein ${\cal F}$-Paar ist, gibt es Schnitte $g'\in{\cal F}(V)$ und $g''\in{\cal F}(U)$ mit $$h=(g')^{-1}g''$$ auf $V\cap U=W$. Schlie\3lich setzen wir $$g_k:=g'h_k\ \:\text{und}\ \:g_l:=g''h_l$$ für $k\in K$ beziehungsweise $l\in I_m$ und erhalten
\begin{eqnarray*}
\text{1.}\ \ \ \:g_i^{-1}g_k= & h_i^{-1}(g')^{-1}g'h_k=h_i^{-1}h_k & =f_{ik}\:,\\
\text{2.}\ \ \ \ g_j^{-1}g_l= & h_j^{-1}(g'')^{-1}g''h_l=h_j^{-1}h_l & =f_{jl}\ \:\text{und}\\
\text{3.}\ \ \ \ g_k^{-1}g_l= & h_k^{-1}(g')^{-1}g''h_l=h_k^{-1}hh_l=h_k^{-1}h_kf_{kl}h_l^{-1}h_l & =f_{kl}\:.
\end{eqnarray*}
Folglich gilt $f=g\Box 1$ für $g=(g_n)_{n\in K\cup I_m}\in{\cal C}^0({\cal U},{\cal F})$, und der Satz ist bewiesen.
\end{proof}

\newpage

\section{Zerfällungssätze}
\label{horst}

\subsection{Der stetige Fall}
\label{hor}

\noindent Sei $G$ die (multiplikative) Gruppe der invertierbaren Elemente einer Banachalgebra $A$. Wir werden mit $G_1$ diejenige Zusammenhangskomponente der Gruppe $G$ bezeichnen, welche das Einselement $1$ enthält. $G_1$ besitzt sogar eine Gruppenstruktur, was man wie folgt einsieht:  

\noindent Für $a\in G_1$ ist $a^{-1}G_1$ wieder zusammenhängend und enthält das Einselement. Also gilt $G_1^{-1}G_1\subset G_1$, und somit ist $G_1$ eine Untergruppe von $G$.

\noindent Weiterhin sei $X$ ein metrischer Raum mit der Abstandsfunktion $\de$, der sich als Vereinigung abzählbar vieler kompakter Teilmengen darstellen lä\3t. $X$ hei\3t {\sl kontrahierbar}\index{kontrahierbare Menge}, wenn eine stetige Abbildung $$H:X\times[0,1]\rightarrow X\ \:\text{mit}\ \:H(x,0)=x\ \:\text{und}\ \:H(x,1)=x_0\in X$$ für alle $x\in X$ existiert. In diesem Falle ist $X$ auch einfach zusammenhängend. 

\noindent Es sei $C^0(X,A)$ die Menge aller auf $X$ stetigen Abbildungen mit Werten in $A$, die offenbar ein Vektorraum über $\mathbb{C}$ ist. Zur Definition einer Topologie auf $C^0(X,A)$ betrachten wir eine Ausschöpfung von $X$ durch kompakte Mengen $K_n$, $n\in\mathbb{N}$ mit $K_n\subset\stackrel{\circ}{K}_{n+1}$ und $X=\bigcup\limits^{\infty}_{n=1}K_n$. Für jedes $n\in\mathbb{N}$ ist durch $$p_n(f):=\|f|_{K_n}\|_{K_n}:=\max\limits_{x\in K_n}\|f(x)\|\:,\ f\in C^0(X,A)$$ eine Halbnorm auf $C^0(X,A)$ erklärt. Durch die Folge $(p_n)_{n\in\mathbb{N}}$ von Halbnormen wird $C^0(X,A)$ zu einem {\sl lokalkonvexen}\index{lokalkonvexer Raum} Raum, dessen Topologie durch die Umgebungen $$W_{n,r}=\{f\in C^0(X,A)\mid p_n(f)<r\}$$ mit $r>0$ und $n\in\mathbb{N}$ der Nullabbildung $0\in C^0(X,A)$ definiert ist. Diese Topologie ist hausdorffsch, da für jede von der Nullabbildung verschiedene Abbildung $f\in C^0(X,A)$ ein $n\in\mathbb{N}$ existiert mit $p_n(f)>0$. Da nur abzählbar viele Halbnormen auf $C^0(X,A)$ betrachtet werden und die lokalkonvexe Topologie hausdorffsch ist, wird der Raum $C^0(X,A)$ metrisierbar. Die Konvergenz bezüglich der oben definierten Topologie ist gerade die gleichmä\3ige Konvergenz von Abbildungen auf allen kompakten Mengen $K_n$, $n\in\mathbb{N}$ und folglich auch auf allen kompakten Teilmengen von $X$. Damit ist die Grenzabbildung beliebiger Cauchyfolgen in $C^0(X,A)$ wieder stetig auf $X$ und die erklärte lokalkonvexe Topologie ist vollständig. Somit ist $C^0(X,A)$ ein {\sl Fréchetraum}\index{Fréchetraum}.

\noindent Ist $Y\subset X$ eine beliebige nichtleere Teilmenge von $X$, so betrachte man auf dem Raum $C^0(Y,A)$ die entsprechende lokalkonvexe Topologie. Ist speziell $Y=K$ kompakt, so lä\3t sich auf $C^0(K,A)$ die Supremumsnorm $$\|f\|_K:=\max\limits_{x\in K}\|f(x)\|$$ erklären, wobei auf einer kompakten Menge das Supremum als Maximum angenommen wird. Der dadurch entstandene Banachraum $C^0(K,A)$ besitzt die Topologie der gleichmä\3igen Konvergenz. Da $C^0(Y,G)$ wieder für beliebiges $Y\subset X$ eine Teilmenge von $C^0(Y,A)$ ist, betrachten wir auf $C^0(Y,G)$ die induzierte Topologie. Dadurch wird $C^0(Y,G)$ zu einer topologischen (multiplikativen) Gruppe. Man bezeichnet mit $C^0_1(Y,G)$ diejenige Zusammenhangskomponente von $C^0(Y,G)$, welche die Einsabbildung $1$ enthält. Offenbar ist $C^0_1(Y,G)$ bezüglich der punktweise definierten Multiplikation eine topologische Gruppe (analog zur Betrachtung von $G$ und $G_1$ zu Beginn des Abschnitts). Weiterhin ist $C^0(Y,G_1)$ eine topologische (multiplikative) Gruppe.

\noindent Bevor wir uns nun dem Faktorisierungssatz (vergleiche Satz \ref{horst-betty}) in seiner vollen Allgemeinheit zuwenden, betrachten wir folgenden Spezialfall: Sei dazu $A:=\mathbb{R}$, d.h. $G:=\mathbb{R}\setminus\{0\}$.

\begin{stz}\label{horst-steffi} Seien $U_1$, $U_2$ offene Teilmengen von $X$ mit einem nichtleeren kontrahierbaren Durchschnitt. Dann gibt es zu jeder stetigen Funktion $f:U_1\cap U_2\rightarrow\mathbb{R}\setminus\{0\}$ stetige Funktionen $f_1:U_1\rightarrow\mathbb{R}\setminus\{0\}$ und $f_2:U_2\rightarrow\mathbb{R}\setminus\{0\}$ mit $f=f_1^{-1}f_2$ auf $U_1\cap U_2$.
\end{stz}

\begin{proof}[{\bf Beweis}] $(U_1,U_2)$ ist eine offene Überdeckung von $U_1\cup U_2$, und es gibt eine {\sl stetige Zerlegung der Eins}\index{stetige Zerlegung der Eins} $(\varphi_1,\varphi_2)$ mit $\supp(\varphi_2)\subset U_1$, $\supp(\varphi_1)\subset U_2$ und $\varphi_1+\varphi_2=1$, $\varphi_i:U_1\cup U_2\rightarrow[0,1]$, $i=1,2$ stetig.

\begin{figure}[htb]
\begin{center}
  \includegraphics*[width=8cm]{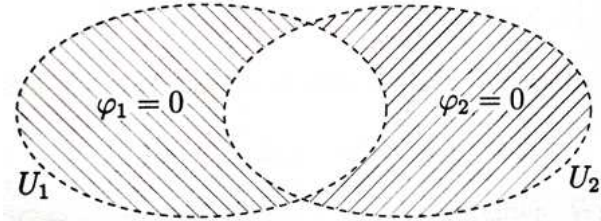}
 \end{center} 
 \caption{\label{fig:bild1} stetige Zerlegung der Eins für $(U_1,U_2)$}
\end{figure}

\noindent Gegeben sei eine stetige Funktion $f:U_1\cap U_2\rightarrow\mathbb{R}\setminus\{0\}$. Wegen der Stetigkeit von $f$ gehören alle $f(x)$ zur gleichen Zusammenhangskomponente $\mathbb{R}_+$ oder $\mathbb{R}_-$. Ohne Beschränkung der Allgemeinheit können wir annehmen, dass die konstante Funktion 1 zu dieser Komponente gehört. Wir betrachten also $\mathbb{R}_+$. Dann kann für jedes $x\in U_1\cap U_2$ das Element $f(x)$ geschrieben werden als $$f(x)=\exp g(x)\:,$$ wobei $g(x)\in\mathbb{R}$ gilt. Die dadurch definierte Funktion $g:U_1\cap U_2\rightarrow\mathbb{R}$ ist stetig. Mit der Zerlegung der Eins gilt: $$f=\exp((\varphi_1+\varphi_2)g)=\exp(\varphi_1g+\varphi_2g)=\exp(\varphi_1g)\exp(\varphi_2g)\:,$$ da $h_1:=\varphi_1g$ und $h_2:=\varphi_2g$ vertauschbar sind, d.h. $h_1h_2=h_2h_1$ gilt. Durch $$g_1:=\begin{cases}\varphi_1g & \text{auf}\ \:U_1\cap U_2 \\ 0 & \text{auf}\ \:U_1\setminus U_2\  \end{cases}\ \ \text{und}\ \ g_2:=\begin{cases}\varphi_2g & \text{auf}\ \:U_1\cap U_2 \\ 0 & \text{auf}\ \:U_2\setminus U_1\  \end{cases}$$ wird $\varphi_1g$ und $\varphi_2g$ auf $U_1$ bzw. $U_2$ stetig fortgesetzt. Setzt man $$f_1:=\exp (-g_1)\ \:\text{und}\ \:f_2:=\exp g_2\:,$$ auf $U_1$ beziehungsweise auf $U_2$, so erhält man $f=f_1^{-1}f_2$ auf $U_1\cap U_2$, und der Satz ist bewiesen.
\end{proof}

\noindent Im allgemeinen Fall unserer betrachteten Gruppe $G$ kann allerdings keine Kommutativität vorausgesetzt werden, und eine neue Beweisidee ist erforderlich. 

\begin{lmm}\label{horst-helena} Seien $\Omega\subset X$ kompakt, $\Omega'\subset$ $\stackrel{\circ}{\Omega}$ abgeschlossen und $f\in C^0_1(\Omega,G)$. Dann gibt es eine Abbildung $\tilde{f}\in C^0_1(X,G)$, so dass $\tilde{f}$ mit $f$ auf $\Omega'$ übereinstimmt.
\end{lmm}

\begin{proof}[{\bf Beweis}] 

\noindent
\begin{enumerate}
\item Es ist bekannt (siehe zum Beispiel \cite{que}, Seite 221), dass eine zusammenhängende topologische Gruppe durch jede offene Umgebung des Einselementes erzeugt werden kann. Speziell betrachten wir die Umgebung $$Z:=\{(1+h)\in C^0_1(\Omega,G)\mid\|h\|_{\Omega}<1\}\:.$$ Jedes Element $f\in C^0_1(\Omega,G)$ lä\3t sich dann in der Form 
\begin{equation}\label{schlaumeise}
f=(1+h_1)\ldots(1+h_m),\ m\geq 1 
\end{equation}
für endlich viele $(1+h_j)\in Z$, $j=1,\ldots,m$ schreiben. Wir wollen an dieser Stelle dennoch den Beweis dafür angeben:

\noindent Es genügt zu zeigen, dass die Menge ${\cal Q}$ aller Abbildungen aus $C^0_1(\Omega,G)$ der Form (\ref{schlaumeise}) offen und abgeschlossen in $C^0_1(\Omega,G)$ ist, und nichtleer. Es sei $g_1\in{\cal Q}$ und die Abbildung $g_2\in C^0_1(\Omega,G)$ erfülle die Bedingung $\|g_2-g_1\|_{\Omega}<\frac{1}{\|g_1^{-1}\|_{\Omega}}$. Für $h:=g_2g_1^{-1}-1$ ergibt sich $g_2=(1+h)g_1$, wobei $$\|h\|_{\Omega}=\|g_2g_1^{-1}-1\|_{\Omega}=\|(g_2-g_1)g_1^{-1}\|_{\Omega}\leq\|g_2-g_1\|_{\Omega}\cdot\|g_1^{-1}\|_{\Omega}<1\:.$$ Folglich besitzt auch $g_2$ die Form (\ref{schlaumeise}), d.h. $g_2\in{\cal Q}$. Also ist die Menge ${\cal Q}$ offen. Eine Folge von Abbildungen $g_n\in{\cal Q}$ konvergiere bezüglich der Norm $\|.\|_{\Omega}$ gegen eine Abbildung $g\in\ C^0_1(\Omega,G)$. Dann ist für hinreichend gro\3es $n$ die Ungleichung $\|g-g_n\|_{\Omega}<\frac{1}{\|g^{-1}_n\|_{\Omega}}$ erfüllt. Wir setzen $h:=gg_n^{-1}-1$. Also gilt $g=(1+h)g_n$, wobei $\|h\|_{\Omega}=\|gg^{-1}_n-1\|_{\Omega}<1$, d.h. $g\in{\cal Q}$, und die Menge ${\cal Q}$ ist abgeschlossen. Wegen $g=1\in{\cal Q}$ ist schlie\3lich ${\cal Q}$ nichtleer.

\newpage

\item Das Lemma von Urysohn liefert uns eine stetige Funktion $$\chi:X\rightarrow[0,1]\ \:\text{mit}\ \:\chi(\Omega')=\{1\}\ \:\text{und}\ \:\chi(\overline{X\setminus\Omega})=\{0\}\:.$$ 

\begin{figure}[htb]
\begin{center}
  \includegraphics*[width=8cm]{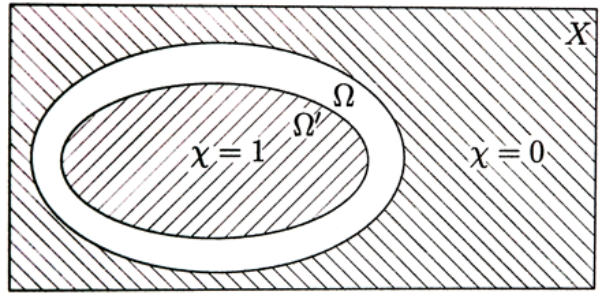}
 \end{center} 
 \caption{\label{fig:bild2} stetige Funktion $\chi$ nach Lemma von Urysohn}
\end{figure}

\noindent Jetzt setzen wir $$\tilde{f}:=\begin{cases}(1+\chi h_1)\ldots(1+\chi h_m) & \text{auf}\ \:\Omega \\ 1 & \text{auf}\ \:X\setminus\Omega \end{cases}\ \ .$$ Dann ist $\tilde{f}$ stetig auf $X$. Weiter sind $\tilde{f}(x)$ für alle $x\in X$ invertierbar, also gilt $\tilde{f}\in C^0(X,G)$, und man erhält $\tilde{f}=f$ auf $\Omega'$.

\item Es bleibt zu zeigen, dass $\tilde{f}\in C^0_1(X,G)$ gilt. Wie bei 1. lä\3t sich $f$ in der Form $f=(1+h_1)\ldots(1+h_m)$ schreiben. Für jedes $j\in\{1,\ldots,m\}$ seien $$\tilde{h}_j:=\begin{cases}\chi h_j & \text{auf}\ \:\Omega \\ 0 & \text{auf}\ \:X\setminus\Omega \end{cases}\ \ .$$ Dann ist die Abbildung $(1+\tilde{h}_j)\in C^0(X,G)$ und $\|\tilde{h}_j(x)\|\leq\|h_j\|_{\Omega}<1$ für alle $x\in X$ und für jedes $j\in\{1,\ldots,m\}$. Es gilt: $$\tilde{f}=(1+\tilde{h}_1)\ldots(1+\tilde{h}_m)\in C^0(X,G)\:.$$ Schlie\3lich ist durch $$\tilde{f}_t:=(1+t\tilde{h}_1)\ldots(1+t\tilde{h}_m)\:,\ t\in [0,1]$$ eine stetige Schar von Abbildungen $\tilde{f}_t\in C^0(X,G)$ erklärt mit $\tilde{f}_1=\tilde f$ und $\tilde{f}_0=1$, so dass wie behauptet $\tilde{f}\in C^0_1(X,G)$ gilt.
\end{enumerate}
\end{proof}

\noindent Das nächste Lemma (vergleiche Lemma 4) ist eine Verallgemeinerung von Lemma \ref{horst-helena}. Dazu benutzen wir folgenden Satz (siehe \cite{lem}, Seite 49 oder \cite{gro}, Seite 12), der auch für einen allgemeinen metrischen Raum $X$ und einen allgemeinen normierten Raum statt der von uns betrachteten Banachalgebra $A$ gilt, was hier allerdings irrelevant ist.

\begin{stz}[{\bf Fortsetzungssatz von Dugundji}]\label{horst-maria}\index{Fortsetzungssatz von Dugundji} Ist $\Omega\neq\emptyset$ eine abgeschlossene Teilmenge von $X$ und $h:\Omega\rightarrow\ A$ stetig, so existiert eine stetige Abbildung $\hat{h}:X\rightarrow A$ mit $\hat{h}|_{\Omega}=h$, wobei $\hat{h}(X)$ in der konvexen Hülle $\mathrm{con}(h(\Omega))$ von $h(\Omega)$ liegt.
\end{stz}

\begin{proof}[{\bf Beweis}] Für jedes $x\in X\setminus\Omega$ sei $$U(x):=K_{\varrho_0}(x)\:,\ \varrho_0:=\frac{1}{4}\dist(x,\Omega)$$ die offene Kugel um $x$ mit dem Radius $\frac{1}{4}\dist(x,\Omega)$. Dann gilt für den Durchmesser von $U(x)$: 
\begin{equation}\label{horst-maria-schnuffi}
\diam U(x)\leq\frac{1}{2}\dist(x,\Omega)\leq\frac{1}{2}\cdot\frac{4}{3}\dist(U(x),\Omega)<\dist(U(x),\Omega)\:.
\end{equation}

\psfrag{distx}{$\dist(x,\Omega)$}
\psfrag{distU}{$\dist(U(x),\Omega)$}
\psfrag{rho0}{$\rho_0$}
\begin{figure}[htb]
\begin{center}
 \psfrag{X}{$X\setminus\Omega$}
 \psfrag{Omega}{$\Omega$}
 \psfrag{U(x)}{$U(x)$}
  \includegraphics*[width=8cm]{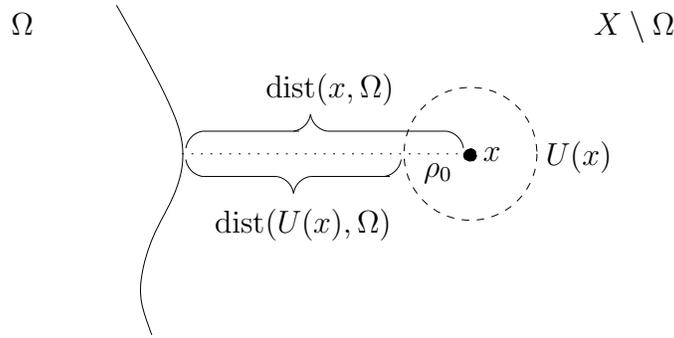}
 \end{center} 
 \caption{\label{fig:bild3} Abstandsbetrachtungen für $U(x)$ und $\Omega$}
\end{figure}

\noindent Wegen $X\setminus\Omega=\bigcup\limits_{x\in X\setminus\Omega}U(x)$ gibt es nach dem Satz von Stone (siehe \cite{lem}, Seite 47) eine {\sl lokal endliche Verfeinerung}\index{Verfeinerung, lokal endliche} der $U(x)$, die ebenfalls $X\setminus\Omega$ überdeckt, etwa durch offene Mengen $V_{\lambda}$, $\lambda\in\Lambda$. Das hei\3t, zu jedem $V_{\lambda}$ existiert ein $U(x)$ mit $V_{\lambda}\subset U(x)$, und jedes $x\in X\setminus\Omega$ besitzt eine offene Umgebung, welche nur mit endlich vielen der $V_{\lambda}$ einen nichtleeren Schnitt hat. Man definiere für $x\in X\setminus\Omega$ die Funktion $$\alpha(x):=\sum_{\lambda\in\Lambda}\dist(x,X\setminus V_{\lambda})\:.$$ Dabei besteht die Summe in einer Umgebung von $x$ aus der gleichen endlichen Anzahl nicht verschwindender Summanden. Somit ist $\alpha$ als Abstandsfunktion stetig und grö\3er als Null in $X\setminus\Omega$. Desweiteren wählen wir zu jedem $\lambda\in\Lambda$ ein festes $\omega_{\lambda}\in\Omega$ mit 
\begin{equation}\label{horst-maria-kleffi}
\dist(V_{\lambda},\omega_{\lambda})<2\dist(V_{\lambda},\Omega).	
\end{equation}
Dies ist wegen $\dist(V_{\lambda},\Omega)\geq\dist(U(x),\Omega)>0$ für $V_{\lambda}\subset U(x)$, $x\in X\setminus\Omega$ möglich. 

\psfrag{omegalambda}{$\omega_{\lambda}$}
\psfrag{Vlambda}{$V_{\lambda}$}
\psfrag{distVomega}{$\dist(V_{\lambda},\omega_{\lambda})$}
\psfrag{distV}{$\dist(V_{\lambda},\Omega)$}
\begin{figure}[htb]
\begin{center}
 \psfrag{X}{$X\setminus\Omega$}
 \psfrag{Omega}{$\Omega$}
 \psfrag{U(x)}{$U(x)$}
  \includegraphics*[width=8cm]{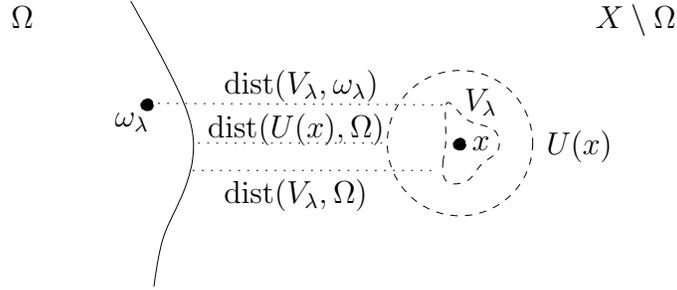}
 \end{center} 
 \caption{\label{fig:bild4} Abstandsbetrachtungen für $V_{\lambda}$ und $\Omega$}
\end{figure}

\noindent Wir setzen: $$\hat{h}(x):=\begin{cases}h(x) & \text{für}\ \:x\in\Omega \\ \frac{1}{\alpha(x)}\sum\limits_{\lambda\in\Lambda}\dist(x,X\setminus V_{\lambda})h(\omega_{\lambda}) & \text{für}\ \:x\in X\setminus\Omega \end{cases}\ \ .$$ Offenbar ist diese Abbildung auf $X$ definiert, und es gilt $\hat{h}(X)\subset\con(h(\Omega))$. Die Abbildung $\hat{h}$ ist sowohl in $\stackrel{\circ}{\Omega}$ als auch in $X\setminus\Omega$ stetig, da $h$ und $\alpha$ stetig sind und die $\omega_{\lambda}$ fest gewählt sind. 

\noindent Es bleibt die Stetigkeit auf dem Rand von $\Omega$ zu zeigen: Seien dazu $x_0\in\partial\Omega$ und $\varepsilon>0$ beliebig gewählt. Zunächst gibt es ein $\delta>0$ mit $$\|\hat{h}(x)-\hat{h}(x_0)\|<\varepsilon\ \text{für}\ \:x\in\Omega\:,\ \de(x,x_0)<\delta\:.$$ Nun sei $x\in X\setminus\Omega$ so gewählt, dass $\de(x,x_0)<\frac{\delta}{4}$ erfüllt ist. Aufgrund der Definition von $\alpha$ und $\hat{h}$, der Gleichheit $\hat{h}(x_0)=\frac{1}{\alpha(x)}\sum\limits_{\lambda\in\Lambda}\dist(x,X\setminus V_{\lambda})\hat{h}(x_0)$ und der Dreiecksungleichung gilt dann die Ungleichung:
\begin{equation}\label{horst-maria-schnucki}
\|\hat{h}(x)-\hat{h}(x_0)\|\leq\frac{1}{\alpha(x)}\sum\limits_{\lambda\in\Lambda}\dist(x,X\setminus V_{\lambda})\|\hat{h}(\omega_{\lambda})-\hat{h}(x_0)\|\:,
\end{equation}
und es genügt folglich zu zeigen, dass für alle $x$, für die $\dist(x,X\setminus V_{\lambda})\neq 0$ und $\de(x,x_0)<\frac{\delta}{4}$ gilt, die Ungleichung $\|\hat{h}(\omega_{\lambda})-\hat{h}(x_0)\|<\varepsilon$ erfüllt ist. Dabei bedeutet $\dist(x,X\setminus V_{\lambda})\neq 0$ aber genau $x\in V_{\lambda}\subset U(x_{\lambda})$, wobei $x_{\lambda}$ passend gewählt sei. Dann gilt:
\begin{eqnarray*}
\de(\omega_{\lambda},x_0) & \leq & \de(\omega_{\lambda},x)+\de(x,x_0)\ \ (\text{wegen der Dreiecksungleichung})\\
 & \leq & \dist(V_{\lambda},\omega_{\lambda})+\diam V_{\lambda}+\de(x,x_0)\ \ (\text{wegen}\ x\in V_{\lambda})\\
 & \leq & 2\dist(V_{\lambda},\Omega)+\diam U(x_{\lambda})+\de(x,x_0)\ \ (\text{wegen (\ref{horst-maria-kleffi}) und } V_{\lambda}\subset U(x_{\lambda}))\\
 & \leq & 2\dist(V_{\lambda},\Omega)+\dist(U(x_{\lambda}),\Omega)+\de(x,x_0)\ \ (\text{wegen (\ref{horst-maria-schnuffi})})\\
 & \leq & 3\dist(V_{\lambda},\Omega)+\de(x,x_0)\ \ (\text{wegen } V_{\lambda}\subset U(x_{\lambda}))\\
 & \leq & 4\de(x,x_0)\ \ (\text{wegen } x\in V_{\lambda }\ \text{und}\ x_0\in\Omega)
\end{eqnarray*}
Somit ist für $\de(x,x_0)<\frac{\delta}{4}$ nun $\de(\omega_{\lambda},x_0)<\delta$, also $\|\hat{h}(\omega_{\lambda})-\hat{h}(x_0)\|<\varepsilon$. Dann ist nach (\ref{horst-maria-schnucki}) auch $\|\hat{h}(x)-\hat{h}(x_0)\|<\varepsilon$.
\end{proof}

\begin{lmm}\label{horst-maja} Seien $\Omega\subset X$ kompakt und $f\in C^0_1(\Omega,G)$. Dann gibt es eine Abbildung $\tilde{f}\in C^0_1(X,G)$, so dass $\tilde{f}$ mit $f$ auf $\Omega$ übereinstimmt.
\end{lmm}

\begin{proof}[{\bf Beweis}] Wie im Beweis von Lemma \ref{horst-helena} bereits gezeigt wurde, lä\3t sich jede Abbildung $f\in C^0_1(\Omega,G)$ in der Form $$f=(1+h_1)\ldots(1+h_m)\:,\ m\geq 1$$ für endlich viele $(1+h_j)\in C^0_1(\Omega,G)$ mit $\|h_j\|_{\Omega}<1$, $j=1,\ldots,m$ schreiben.

\noindent Nach dem eben bewiesenen Fortsetzungssatz von Dugundji sind auf einer kompakten Menge $\Omega\subset X$ die Abbildungen $h_j:\Omega\rightarrow A$ stetig auf $X$ fortsetzbar, d.h. es gibt stetige Abbildungen $\hat{h}_j:X\rightarrow A$ mit $\hat{h}_j|_{\Omega}=h_j.$ Für $\|h_j\|_{\Omega}<1$ gilt auch $\|\hat{h}_j(x)\|\leq\|h_j\|_{\Omega}<1$ für alle $x\in X$ aufgrund der Bedingung für die konvexe Hülle. Wir setzen: $$\tilde{f}:=(1+\hat{h}_1)\ldots(1+\hat{h}_m)\in C^0(X,G)\:.$$ Dann ist analog wie im Beweis von Lemma \ref{horst-helena} durch $$\tilde{f}_t:=(1+t\hat{h}_1)\ldots(1+t\hat{h}_m)\:,\ t\in [0,1]$$ eine stetige Schar von Abbildungen $\tilde{f}_t\in C^0(X,G)$ erklärt mit $\tilde{f}_1=\tilde f$ und $\tilde{f}_0=1$, so dass wie behauptet $\tilde{f}\in C^0_1(X,G)$ gilt.
\end{proof}

\begin{lmm}[{\bf Fortsetzungssatz}]\label{horst-karla} Seien $W\subset X$ abgeschlossen und kontrahierbar sowie $f\in C^0(W,G_1)$. Dann gibt es eine Abbildung $\tilde{f}\in C^0_1(X,G)$, so dass $\tilde{f}$ mit $f$ auf $W$ übereinstimmt.
\end{lmm}

\begin{proof}[{\bf Beweis}] 

\noindent
\begin{enumerate}
\item \noindent Zunächst zeigt man, dass aus der Kontrahierbarkeit von $W$ und dem Zusammenhang von $G_1$ der Zusammenhang von $C^0(W,G_1)$ folgt:

\noindent Sei $f$ eine beliebige Abbildung aus $C^0(W,G_1)$. Da $W$ kontrahierbar ist, existiert eine stetige Abbildung $$H:W\times [0,1]\rightarrow W\ \:\text{mit}\ \:H(x,0)=x\ \:\text{und}\ \:H(x,1)=x_0\in W$$ für alle $x\in W$. Man definiert eine Abbildungsschar $$f_t(x):=f(H(x,t))\:,\ x\in W\:,\ t\in [0,1]\:,$$ wobei $H$ die eben erklärte Kontraktionsabbildung ist. Die Abbildung $t\mapsto f_t$ ist stetig, und verbindet in $C^0(W,G_1)$ das gegebene Element $f_0=f$ mit dem Element $f_1$. Da die Abbildung $f_1\equiv f(x_0)$ konstant ist, und $f(x_0)\in G_1$ gilt, folgt der Zusammenhang von $C^0(W,G_1)$. Also gilt $C^0(W,G_1)\subset C^0_1(W,G)$, und folglich auch $f\in C^0_1(W,G)$. 
\item Man wähle kompakte Mengen $\Omega_n\subset$ $\stackrel{\circ}{\Omega}_{n+1}$ mit $X=\bigcup\limits^{\infty}_{n=1}\Omega_n$ so, dass die Mengen $W\cap\Omega_n$ für alle $n\in\mathbb{N}$ kontrahierbar sind, wobei $\Omega_1$ mit $W$ einen nichtleeren Schnitt hat. Für jedes $n\in\mathbb{N}$ existiert eine stetige Abbildung $H_n:W\cap\Omega_{n+1}\times[0,1]\rightarrow W\cap\Omega_{n+1}$ mit
\begin{itemize}
	\item $H_n(x,0)=x\ \text{und}\ H_n(x,1)=x_0\in W\cap\Omega_n\ \text{für alle}\ x\in W\cap\Omega_{n+1}\ \text{sowie}$
	\item $H_n(W\cap\Omega_n,t)\subset W\cap\Omega_n\ \text{für alle}\ t\in[0,1].$
\end{itemize}
Man konstruiere induktiv eine Folge von Abbildungen $f_n\in C^0_1(X,G)$ mit:
\begin{enumerate}
\item $f_n=f$ auf $W\cap \Omega_n$ für $n\geq 1$ und
\item $f_n=f_{n-1}$ auf $\Omega_{n-1}$ für $n\geq 2$.
\end{enumerate}

\begin{figure}[htb]
\begin{center}
  \includegraphics*[width=12cm]{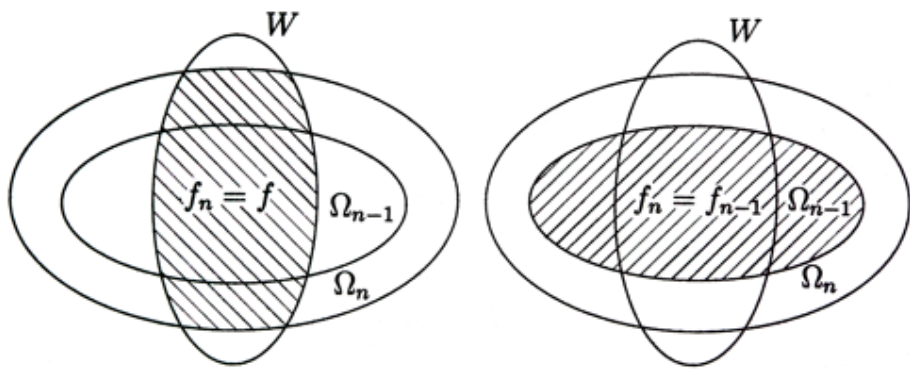}
 \end{center} 
 \caption{\label{fig:bild5} zur induktiven Konstruktion der Abbildungen $f_n$}
\end{figure}

\noindent Für den Induktionsanfang wenden wir Lemma \ref{horst-maja} auf $\Omega:=W\cap\Omega_1$ an und erhalten $f_1$ als stetige Fortsetzung von $f|_{W\cap \Omega_1}$ auf $X$.

\noindent Beim Induktionsschritt seien nun $f_1,\ldots,f_k\in C^0_1(X,G)$, $k\geq 1$ bereits konstruiert und es gelte:
\begin{enumerate}
\item $f_n=f$ auf $W\cap \Omega_n$ für $n=1,\ldots,k$ und
\item $f_n=f_{n-1}$ auf $\Omega_{n-1}$ für $n=2,\ldots,k$.
\end{enumerate}
Die Konstruktion von $f_{k+1}$ erfolgt dann in 2 Schritten. Man wende zunächst Lemma \ref{horst-maja} auf $\Omega:=W\cap\Omega_{k+1}$ an und erhält eine Abbildung $\tilde{f}_{k+1}\in C^0_1(X,G)$ mit $\tilde{f}_{k+1}=f$ auf $W\cap\Omega_{k+1}$. Im zweiten Schritt setze man $$g:=\begin{cases}f_k\tilde{f}_{k+1}^{-1} & \text{auf}\ \:\Omega_k \\ 1 & \text{auf}\ \:W\cap\Omega_{k+1} \end{cases}\ \ \in C^0(\Omega_k\cup(W\cap\Omega_{k+1}),G)\:.$$ Diese Definition ist möglich, da die Abbildung $g$ wegen (a) auf dem gemeinsamen Durchschnitt $W\cap\Omega_k$ gleich $1$ ist. Als nächstes weisen wir nach, dass $g$ aus $C^0_1(\Omega_k\cup(W\cap\Omega_{k+1}),G)$ ist:

\noindent Dazu setzen wir $$h:=g\tilde{f}_{k+1}f_k^{-1}\in C^0(\Omega_k\cup(W\cap\Omega_{k+1}),G)$$ und zeigen $h\in C^0(\Omega_k\cup(W\cap\Omega_{k+1}),G)$. Es gilt: $$h=\begin{cases}1 & \text{auf}\ \:\Omega_k \\ \tilde{f}_{k+1}f_k^{-1} & \text{auf}\ \:W\cap\Omega_{k+1} \end{cases}\ \ .$$ Mit der Kontraktionsabbildung $H_k$ definieren wir die stetige Schar $$h_t:=\begin{cases}1 & \text{auf}\ \:\Omega_k \\ h(H_k(.,t)) & \text{auf}\ \:W\cap\Omega_{k+1} \end{cases}\ \ \in C^0(\Omega_k\cup(W\cap\Omega_{k+1}),G)\:,\ t\in[0,1]\:.$$ Diese Definition ist möglich, da auf der gemeinsamen Schnittmenge $W\cap\Omega_k$ die Inklusion $H_k(W\cap\Omega_k,t)\subset W\cap\Omega_k\subset\Omega_k$ gilt und folglich $h(H_k(x,t))=1$ für alle $x\in W\cap\Omega_k$ und alle $t\in [0,1]$ laut Definition von $h$ erfüllt ist. Weiterhin gilt mit Hilfe der Eigenschaften von $H_k$: $$h_0=\begin{cases}1 & \text{auf}\ \:\Omega_k \\ h & \text{auf}\ \:W\cap\Omega_{k+1} \end{cases}\ =h\ \ \text{und}\ \ h_1=1\:,$$ woraus wie behauptet $h$, und damit auch $g\in C^0_1(\Omega_k\cup(W\cap\Omega_{k+1}),G)$ folgt.

\noindent Aus Lemma \ref{horst-maja} folgt nun, dass es eine Abbildung $\tilde{g}\in C^0_1(X,G)$ gibt mit $\tilde{g}=g$ auf $\Omega_k\cup(W\cap\Omega_{k+1})$. Durch $$f_{k+1}:=\tilde{g}\tilde{f}_{k+1}$$ ist sowohl $f_{k+1}=f$ auf $W\cap\Omega_{k+1}$, als auch $f_{k+1}=f_k$ auf $\Omega_k$ erfüllt.
\item Die Folge der Abbildungen $(f_n)_{n\in{\mathbb N}}$ konvergiert auf $X$ punktweise gegen eine Abbildung $\tilde{f}$ mit $\tilde{f}(x)=f_n(x)$ für $x\in\Omega_n$ nach (b). Da die Abbildungen $f_n$ auf $X$ stetig sind und die Konvergenz auf jeder der kompakten Mengen $\Omega_n$ gleichmä\3ig ist, ist die Grenzabbildung $\tilde{f}$ auf allen $\Omega_n$ und damit auch auf $X$ stetig. Folglich gilt $\tilde{f}\in C^0_1(X,G)$ und $\tilde{f}=f$ auf $W$ nach (a). 
\end{enumerate}
\end{proof}

\noindent Nach den eben betrachteten Lemmata kann man jetzt den nächsten Satz beweisen:

\begin{stz}[{\bf Faktorisierungssatz}]\label{horst-betty} Seien $U_1$ und $U_2$ offene Teilmengen von $X$ mit einem nichtleeren kontrahierbaren Durchschnitt. Dann gibt es zu jeder stetigen Abbildung $f:U_1\cap U_2\rightarrow G$ stetige Abbildungen $f_1:U_1\rightarrow G$ und $f_2:U_2\rightarrow G$ mit $f=f_1^{-1}f_2$ auf $U_1\cap U_2$.
\end{stz}

\begin{bmr} Man sagt auch, $f$ lä\3t sich {\sl faktorisieren}\index{Faktorisierung von stetigen Abbildungen}.
\end{bmr}

\begin{proof}[{\bf Beweis}] Man wähle $W_1\subset U_1$ und $W_2\subset U_2$ so, dass sie abgeschlossen in $U_1\cup U_2$ sind, ihr Durchschnitt wieder nichtleer und kontrahierbar ist, sowie $W_1\cup W_2=U_1\cup U_2$ gilt. Dann ist $W_1\cap W_2$ abgeschlossen in $U_1\cup U_2$, und es gilt $W_1\cap W_2\subset U_1\cap U_2$.
 
\psfrag{W1}{$W_1$}
\psfrag{W2}{$W_2$}
\psfrag{U1}{$U_1$}
\psfrag{U2}{$U_2$} 
\begin{figure}[htb]
\begin{center}
  \includegraphics*[width=7cm]{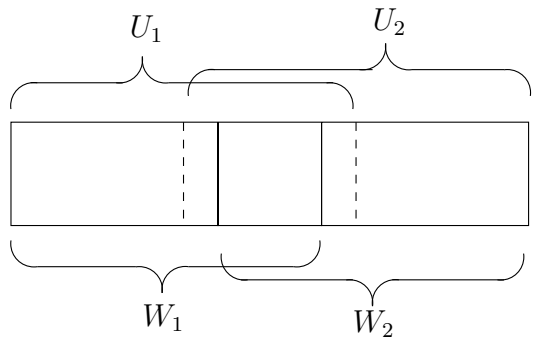}
 \end{center} 
 \caption{\label{fig:bild6} relativ abgeschlossene Teilmengen $W_1$ und $W_2$ in $U_1\cup U_2$}
\end{figure} 
 
\noindent Dass eine solche Wahl möglich ist, zeigt folgende Betrachtung: Nach dem Lemma von Urysohn gibt es eine stetige Funktion $$\chi:U_1\cup U_2\rightarrow[0,1]\ \:\text{mit}\ \:\chi(U_1\setminus U_2):=\{0\}\ \:\text{und}\ \:\chi(U_2\setminus U_1):=\{1\}\:.$$ Nun seien $$W_1:=\chi^{-1}\Big(\Big[0,\frac{3}{4}\Big]\Big)\ \:\text{und}\ \:W_2:=\chi^{-1}\Big(\Big[\frac{1}{4},1\Big]\Big)\:.$$ Als Urbilder abgeschlossener Mengen bezüglich einer stetigen Funktion sind $W_1$ und $W_2$ abgeschlossen. Weiterhin gilt $W_1\cup W_2=\chi^{-1}([0,1])=U_1\cup U_2$ und $W_1\cap W_2=\chi^{-1}([\frac{1}{4},\frac{3}{4}])\neq\emptyset$. Ohne Beschränkung der Allgemeinheit kann man annehmen, dass $W_1\cap W_2$ kontrahierbar ist. Andernfalls betrachte man eine Zusammenhangskomponente $V$ von $(U_1\cap U_2)\setminus(W_1\cap W_2)$. Gilt $\chi(\partial V)=\{c\}$, $c\in\{\frac{1}{4},\frac{3}{4}\}$, so definieren wir $\chi$  durch $\chi(V)=\{c\}$ stetig um. 

\noindent Au\3erdem können wir ohne Beschränkung der Allgemeinheit $f(x_0)=1$ für einen beliebigen aber festen Punkt $x_0\in U_1\cap U_2$ annehmen. Andernfalls betrachte man statt $f$ die Abbildung $(f(x_0))^{-1}f(x)$, $x\in U_1\cap U_2$. 

\noindent Dann gehören alle Werte $f(x)$ für $x\in U_1\cap U_2$ zu $G_1$, also $f\in C^0(U_1\cap U_2,G_1)$. Für die Anwendung des Fortsetzungssatzes (Lemma \ref{horst-karla}) seien $W:=W_1\cap W_2$ und $X:=U_1\cup U_2$. Nun folgt, dass es eine stetige Abbildung $\tilde{f}:U_1\cup U_2\rightarrow G$ gibt mit $\tilde{f}=f$ auf $W_1\cap W_2$. Man setze: $$\tilde{f}_1:=1\ \:\text{und}\ \:\tilde{f}_2:=\tilde{f}$$ auf $U_1\cup U_2$ und es gilt $f=\tilde{f}_1^{-1}\tilde{f}_2$ auf $W_1\cap W_2$. Es bleibt zu zeigen, dass $f$ auch auf $U_1\cap U_2$ faktorisiert werden kann. Dazu setze man $$f_2:=\begin{cases}f & \text{auf}\ \:W_1\cap U_2 \\ \tilde{f} & \text{auf}\ \:W_2 \end{cases}\ \ \text{und}\ \ f_1:=\begin{cases}f_2f^{-1} & \text{auf}\ \:U_1\cap U_2 \\ 1 & \text{auf}\ \:U_1\setminus W_2 \end{cases}\ \ .$$ 

\begin{figure}[htb]
\hspace{4cm}
   \includegraphics*[width=12cm]{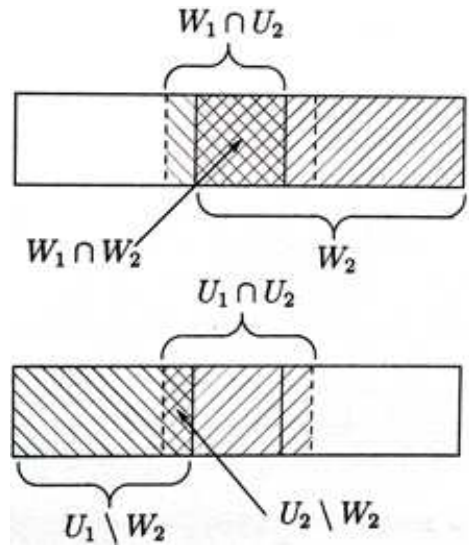}
 \caption{\label{fig:bild7} Schnittmengendarstellung zur Definition von $f_2$ beziehungsweise $f_1$}
\end{figure} 
 
\noindent Dies ist möglich, da im ersten Fall $f=\tilde{f}$ auf der gemeinsamen Schnittmenge $W_1\cap W_2$, im zweiten Fall $f_2=f$ und folglich $f_1=1$ auf der gemeinsamen Schnittmenge $U_2\setminus W_2$ gilt. Die Abbildungen sind wieder stetig und es sind $f_1(z)$ für alle $z\in U_1$ und $f_2(z)$ für alle $z\in U_2$ invertierbar. Schlie\3lich gilt auf $U_1\cap U_2$ : $$f_1^{-1}f_2=(f_2f^{-1})^{-1}f_2=ff_2^{-1}f_2=f\:.$$
\end{proof}

\noindent Es folgt leicht:

\begin{stz}[{\bf Zerfällungssatz}]\label{horst-petra} Es sei ${\cal C}^G$ die Garbe der Keime stetiger Abbildungen mit Werten in $G$ und $U_1,\ldots,U_n$, $n\in{\mathbb N}$ seien offene Mengen in $X$, so dass die Durchschnitte $(U_1\cup\ldots\cup U_k)\cap U_{k+1}$, $k=1,\ldots,n-1$ kontrahierbar sind. 

\noindent Dann gilt $H^1({\cal U},{\cal C}^G)=1$ für die Überdeckung ${\cal U}:=(U_k)_{1\leq k\leq n}$ von $U_1\cup\ldots\cup U_n$.
\end{stz}

\begin{proof}[{\bf Beweis}] Es sei ${\cal F}={\cal C}^G$. Für $n=1$ ist die Aussage trivial. Für $n>1$ ist $(U_1,\ldots,U_n)$ eine $\cal F$-Kette, da $(U_1\cup\ldots\cup U_k,U_{k+1})$ für jedes $k\in\{1,\ldots,n-1\}$  ein $\cal F$-Paar ist. Genau das haben wir in Satz \ref{horst-betty} gezeigt. Nach Satz \ref{harry-kathi} ist folglich die zugehörige Kohomologiemenge trivial, d.h. es gilt $H^1({\cal U},{\cal C}^G)=1$, was zu beweisen war.
\end{proof}

\subsection{Der holomorphe Fall}
\label{st}

\noindent Es sei $G$ wieder die (multiplikative) Gruppe der invertierbaren Elemente einer Banachalgebra $A$ und $G_1$ diejenige Zusammenhangskomponente der Gruppe $G$, welche das Einselement $1$ enthält. Wir wählen nun $X:=\mathbb{C}$ und nennen nichtleere offene Mengen der Form $$R:=\{z\in{\mathbb C}\mid a<\re z<b,c<\im z<d\}\ \:\text{mit}\ \:a,b,c,d\in{\mathbb R}$$ {\sl Rechtecke}\index{Rechteck in $\mathbb C$} in $\mathbb{C}$. Für ein beliebiges Rechteck $R\subset\mathbb{C}$ sei $\overline{\cal O}(R,A)$ der Banachraum aller auf $\overline R$ stetigen und in $R$ holomorphen Abbildungen mit Werten in der Banachalgebra $A$. Als Norm verwendet man die Supremumsnorm $\|.\|_{\overline R}$, da die betrachteten Abbildungen auf der kompakten Menge $\overline R$ stetig sind. Wir weisen darauf hin, dass die Abbildungen aus dem Raum $\overline{\cal O}(R,G)$ auch noch invertierbare Werte auf dem Rand von $R$ besitzen. Weiterhin sei $C^{\infty}(R,A)$ der Raum aller in $R$ beliebig oft reell differenzierbaren Abbildungen mit Werten in $A$.

\begin{lmm}\label{hector-elisa} Seien $R_1$ und $R_2$ zwei Rechtecke in $\mathbb C$, so dass ihre Vereinigung $R_1\cup R_2$ auch wieder ein Rechteck ist. Dann lä\3t sich jede Abbildung $f\in\overline{\cal O}(R_1\cap R_2,A)$ darstellen als $f=f_1-f_2$, wobei $f_1\in\overline{\cal O}(R_1,A)$ und $f_2\in\overline{\cal O}(R_2,A)$ sind. 
\end{lmm}

\begin{bmr} Für zwei Rechtecke, deren Vereinigung wieder ein Rechteck ist, ist auch deren Durchschnitt ein Rechteck.
\end{bmr}

\begin{proof}[{\bf Beweis von Lemma \ref{hector-elisa}}] Es sei $$\chi:\overline{R_1\cup R_2}\rightarrow [0,1]\ \:\text{mit}\ \:\chi(\overline{R_1\setminus R_2})=\{0\}\ \:\text{und}\ \:\chi(\overline{R_2\setminus R_1})=\{1\}$$ eine $C^{\infty}$-Funktion (für die Zerlegung der Eins siehe \cite{tri}, Seite 47). Man setze die Abbildungen $$g_1:=\begin{cases}\chi f & \text{auf}\ \:\overline{R_1\cap R_2} \\ 0 & \text{auf}\ \:\overline R_1\setminus R_2 \end{cases}\ \ ,$$ $$g_2:=\begin{cases}(\chi-1)f & \text{auf}\ \:\overline{R_1\cap R_2} \\ 0 & \text{auf}\ \:\overline R_2\setminus R_1 \end{cases}\ \ ,$$ womit $g_j\in C^{\infty}(R_j,A)$ gilt und $g_j$ auf $\overline R_j,\ j=1,2$ stetig ist. Man erhält $f=g_1-g_2$ auf $\overline{R_1\cap R_2}$. Auf $R_1\cap R_2$ ist das komplexe Vektorfeld $$\frac{\partial}{\partial\overline z}:=\frac{1}{2}\Big(\frac{\partial}{\partial x}+i\frac{\partial}{\partial y}\Big)$$ mit $x:=\re z$ und $y:=\im z$ definiert. Die Abbildung $f$ ist genau dann holomorph auf $R_1\cap R_2$, wenn dort $$\frac{\partial f}{\partial\overline z}\equiv 0$$ gilt, denn die Beziehung $\frac{\partial f}{\partial\overline z}\equiv 0$ besagt, dass die Abbildung $f$ bezüglich der komplexen Variablen $z$ die aus der Funktionentheorie bekannten Cauchy-Riemannschen Differentialgleichungen erfüllt. Aus $f=g_1-g_2$ folgt nun $$\frac{\partial g_1}{\partial\overline z}=\frac{\partial g_2}{\partial\overline z}$$ auf $R_1\cap R_2$, und wir können $$\varphi:=\begin{cases}\frac{\partial g_1}{\partial\overline z} & \text{auf}\ \:R_1 \\ \frac{\partial g_2}{\partial\overline z} & \text{auf}\ \:R_2 \end{cases}\ \ \in C^{\infty}(R_1\cup R_2,A)$$ definieren. Das Pompeiu-Integral (siehe \cite{gun}, Seite 25 oder \cite{wer}, Seite 32) liefert eine Lösung $$h(z):=-\frac{1}{\pi}\int\limits_{R_1\cup R_2}\frac{\varphi(\zeta)}{\zeta-z}\de\lambda(\zeta)\:,\ z\in\overline{R_1\cup R_2}$$ von $\frac{\partial h}{\partial\overline z}=\varphi$, wobei $\de\lambda$ das Lebesque-Ma\3 auf $\mathbb C$ ist. Offensichtlich ist die Abbildung $h\in C^{\infty}(R_1\cup R_2,A)$ und stetig auf $\overline{R_1\cup R_2}$, da der Integrand eine Singularität höchstens vom 1. Grad aufweist. Die auf $\overline R_j$ stetigen Abbildungen $$f_j:=g_j-h\:,\ j=1,2$$ erfüllen $\frac{\partial f_j}{\partial\overline z}\equiv 0$ auf $R_j$ und liefern das Gewünschte.
\end{proof}

\noindent Wir weisen nun nach, dass sich gewisse Abbildungen nahe der Einsabbildung {\sl faktorisieren}\index{Faktorisierung von holomorphen Abbildungen} lassen (siehe Lemma \ref{hector-amata}). Im Beweis benutzen wir den Satz von Graves (siehe \cite{lan}, Seite 193), der hier in leicht vereinfachter Form bewiesen wird.

\begin{stz}[{\bf Graves}]\label{hector-sybile}\index{Satz von Graves} Es sei $V$ eine offene Menge in dem Banachraum $E_1$. Sei $\Theta:V\rightarrow E_2$ eine stetig differenzierbare Abbildung in den Banachraum $E_2$ und $b_0\in V$. Wenn $(D\Theta)(b_0)$ surjektiv ist, dann enthält $\Theta(V)$ eine offene Umgebung von $\Theta(b_0)$. 
\end{stz}

\begin{proof}[{\bf Beweis}] Wir können annehmen, dass $b_0=0$ und $\Theta(b_0)=0$ gilt, da man ansonsten eine Translation in den Banachräumen $E_1$ und $E_2$ durchführen kann. Es genügt zu zeigen, dass für eine offene Kugel $K_r(0)\subset E_1$ um den Mittelpunkt $0$, die Menge $\Theta(K_r(0))\subset E_2$ eine offene Umgebung der $0$ enthält. Sei die lineare Abbildung $$\Lambda:=(D\Theta)(0)$$ die Ableitung von $\Theta$ an der Stelle $b_0=0$. Da $\Lambda$ nach Voraussetzung surjektiv ist, können wir das Open Mapping Theorem (siehe zum Beispiel \cite{lan}, Seite 183) anwenden und erhalten $\Lambda(K_1(0))$ als eine offene Umgebung in $E_2$. Sie enthält wegen der Homogenität von $\Lambda$ die 0, und weiterhin eine abgeschlossene Kugel $\overline K'_{\delta}(0)$ mit $\delta>0$. Folglich gilt $\Lambda(\overline K_1(0))\supset\Lambda(K_1(0))\supset\overline K'_{\delta}(0)=\delta\overline K'_1(0)$. Durch geeignete Wahl der Norm auf $E_2$ können wir $\Lambda(\overline K_1(0))\supset\overline K'_1(0)$ annehmen. Aufgrund der Homogenität von $\Lambda$ gilt auch $\Lambda(\overline K_{\varrho}(0))\supset\overline K'_{\varrho}(0)$ für jedes $\varrho>0$. Für ein beliebig gegebenes $y\in E_2$ mit $\underbrace{\|y\|}_{=:\varrho_0}\leq 1$ existiert nun ein $x\in\overline K_{\varrho_0}(0)\subset E_1$, so dass gilt:
\begin{equation}\label{stern}
\Lambda x=y\ \:\text{und}\ \:\|x\|\leq\varrho_0=\|y\|\:.
\end{equation}
\noindent Es sei $0<\varepsilon<1$ beliebig gewählt. Man nehme $K_r(0)$ mit einem genügend kleinen Radius $r$, so dass uns der Mittelwertsatz (siehe zum Beispiel \cite{lan}, Seite 103) für $x$ und $x'$ aus $\frac{1}{1-\varepsilon}K_r(0)\!:$
\begin{equation}\label{sternchen}
\|\Theta(x)-\Theta(x')-\Lambda(x-x')\|\leq \varepsilon\|x-x'\|
\end{equation}
liefert. Es genügt zu zeigen, dass $$\Theta\Big(\frac{1}{1-\varepsilon}K_r(0)\Big)\supset\Lambda(K_r(0))$$ gilt, was wir im Folgenden tun werden. Sei dazu $x_1\in K_r(0)\subset\frac{1}{1-\varepsilon}K_r(0)$ beliebig gewählt und $y_1:=\Lambda x_1\in\Lambda(K_r(0))$. Wir suchen also ein $x_0\in\frac{1}{1-\varepsilon}K_r(0)$, welches $\Theta(x_0)=y_1$ erfüllt.

\noindent {\bf 1. Schritt}: Wir setzen $y_2:=y_1-\Theta(x_1)$. Wegen (\ref{sternchen}) gilt: $$\|y_1-\Theta(x_1)\|=\|\Lambda x_1-\Theta(x_1)\|\leq \varepsilon\|x_1\|\:.$$ Es existiert ein $x_2$ mit $\Lambda x_2=y_2$ und $\|x_2\|\leq \|y_2\|\leq\varepsilon\|x_1\|$ aufgrund von (\ref{stern}). Dann ist $$x_1+x_2\in(1+\varepsilon)K_r(0)\:,$$ und wegen (\ref{sternchen}) gilt: $$\|y_1-\Theta(x_1+x_2)\|=\|\Theta(x_1)-\Theta(x_1+x_2)+\Lambda x_2\|\leq \varepsilon\|x_2\|\leq \varepsilon^2\|x_1\|\:.$$

\noindent {\bf n-ter Schritt}: Jetzt seien $y_n$ und $x_n$ gegeben mit $\Lambda x_n=y_n$, $\|x_n\|\leq\varepsilon^{n-1}\|x_1\|$ und $$\|y_1-\Theta(x_1+\ldots+x_n)\|\leq\varepsilon^n\|x_1\|\:.$$ Es sei $y_{n+1}:=y_1-\Theta(x_1+\ldots+x_n)$. Wir können nach (\ref{stern}) ein $x_{n+1}$ finden, so dass $\Lambda x_{n+1}=y_{n+1}$ und $\|x_{n+1}\|\leq\|y_{n+1}\|\leq\varepsilon^n\|x_1\|$ erfüllt wird. Dann ist $$x_1+x_2+\ldots+x_{n+1}\in(1+\varepsilon+\ldots+\varepsilon^n)K_r(0)\:,$$ und mit Hilfe von (\ref{sternchen}) folgt wieder
\begin{eqnarray*}
\|y_1 & - & \Theta(x_1+\ldots+x_{n+1})\|\\
 & = & \|y_1-\Theta(x_1+\ldots+x_n)+\Theta(x_1+\ldots+x_n)-\Theta(x_1+\ldots+x_{n+1})\|\\
 & = & \|\Lambda x_{n+1}+\Theta(x_1+\ldots+x_n)-\Theta(x_1+\ldots+x_{n+1})\|\\
 & \leq & \varepsilon\|x_{n+1}\|\leq \varepsilon^{n+1}\|x_1\|\:.
\end{eqnarray*} 
Wir setzen $x_0:=\sum\limits^{\infty}_{k=1}x_k$, wobei diese Reihe wegen $\|x_k\|\leq\varepsilon^{k-1}\|x_1\|$ für alle $k>1$ der Norm nach konvergiert. Man erhält $\Theta(x_0)=y_1$ wegen $\|y_1-\Theta(\sum\limits_{k=1}^{n}x_k)\|\rightarrow 0 $ bei $n\rightarrow\infty$. Weiterhin gilt $x_0\in\frac{1}{1-\varepsilon}K_r(0)$ wegen $\|x_0\|\leq\|x_1\|\sum\limits^{\infty}_{k=1}\varepsilon^{k-1}=\frac{1}{1-\varepsilon}\|x_1\|$, und unser Satz ist bewiesen.
\end{proof}

\begin{lmm}\label{hector-amata} Es gibt ein $\varepsilon_0>0$, so dass für alle Abbildungen $f\in\overline{\cal O}(R_1\cap R_2,G)$ mit $\|f-1\|_{\overline{R_1\cap R_2}}<\varepsilon_0$ Abbildungen $f_1\in\overline{\cal O}(R_1,G)$ und $f_2\in\overline{\cal O}(R_2,G)$ existieren mit $f=f_1^{-1}f_2$ auf $R_1\cap R_2$.
\end{lmm} 

\begin{proof}[{\bf Beweis}] Man setze $B:=\overline{\cal O}(R_1,A)\oplus \overline{\cal O}(R_2,A),\ F:=\overline{\cal O}(R_1\cap R_2,A)$ und definiere eine Abbildung $$\Phi:B\rightarrow F\ \:\text{durch} \ \:\Phi(g_1,g_2):=(1+g_1)(1+g_2)=1+g_1+g_2+g_1g_2\:,$$ wobei $(g_1,g_2)\in B$ ist. Wir erhalten sofort $\Phi(0,0)=1\in\mathrm{Im}\Phi$.

\noindent Seien $b=(b_1,b_2)\in B$ und $h=(h_1,h_2)\in B$. Gesucht ist eine lineare Abbildung $D_b\Phi\in L(B,F)$ mit $\|\Phi(b+h)-\Phi(b)-(D_b\Phi)h\|=o(\|h\|)$ :

\noindent Es ergibt sich $\Phi(b+h)-\Phi(b)=(b_2+1)h_1+(b_1+1)h_2+h_1h_2$. Die durch $$(D_b\Phi)h:=(b_2+1)h_1+(b_1+1)h_2$$ definierte lineare Abbildung erfüllt die Ungleichung: $\|\Phi(b+h)-\Phi(b)-(D_b\Phi)h\|=\|h_1h_2\|\leq\|h_1\|\cdot\|h_2\|\leq\|h\|^2=o(\|h\|),$ so dass $\Phi$ auf $B$ stetig differenzierbar ist, und $D_b\Phi$ die Ableitung von $\Phi$ im Punkt $b$ ist. Speziell erhalten wir für $b=(0,0)$ jetzt $$\Psi(g_1,g_2):=(D_{(0,0)}\Phi)(g_1,g_2)=g_1+g_2\:.$$ Die so erklärte Abbildung $\Psi:\overline{\cal O}(R_1,A)\oplus \overline{\cal O}(R_2,A)\rightarrow \overline{\cal O}(R_1\cap R_2,A)$ ist nach Lemma \ref{hector-elisa} surjektiv. Nun haben wir alle erforderlichen Voraussetzungen, um den Satz von Graves anzuwenden:

\noindent Wir verwenden $E_1:=B,\ E_2:=F,\ V:=\{ (g_1,g_2)\in B\mid\|g_j\|_{\overline R_j}<1,\ j=1,2\},$ $\ \Theta:=\Phi,\ b_0:=(0,0)\ \text{und}\ (D\Theta)(b_0):=\Psi.$ Das liefert uns ein $\varepsilon_0>0$, so dass alle $f\in F$ mit $\|f-1\|_{\overline{R_1\cap R_2}}<\varepsilon_0$ ein Urbild $(g_1,g_2)\in V$ bezüglich $\Phi$ besitzen, wobei $f=\Phi(g_1,g_2)=(1+g_1)(1+g_2)$ gilt. Die Abbildungen $$f^{-1}_1:=1+g_1\ \:\text{und}\ \:f_2:=1+g_2$$ erfüllen das Geforderte\footnote{Dabei sind in einer Banachalgebra die $f_j(z)$ $\forall z\in\overline R_j$, $j=1,2$ durch die Wahl von $V$ invertierbar (siehe \cite{lan}, Seite 67).}, und damit ist das Lemma bewiesen.
\end{proof}

\noindent Die Notwendigkeit dieser zunächst sehr speziell erscheinenden Aussage von Lemma \ref{hector-amata} wird im Beweis von Satz \ref{hector-manu} deutlich. Vorausgeschickt werden drei Hilfssätze und folgende Redeweise: 

\noindent Sei $Y$ eine beliebige nichtleere Menge in $\mathbb C$. Eine Abbildung $f$ hei\3t {\sl holomorph in}\index{holomorph in} $Y$, wenn $f$ in einer offenen Menge $U$ mit $U\supset Y$ definiert und holomorph ist (Definition vergleiche \cite{gra}, Seite 48). Falls $Y$ selbst offen ist, so kann natürlich $U=Y$ gewählt werden.

\noindent Wir weisen darauf hin, dass für jedes Rechteck $R\subset{\mathbb C}$ und jede holomorphe Abbildung $f:\overline R\rightarrow G$ auch $f\in\overline{\cal O}(R,G)$ gilt, die Umkehrung aber im Allgemeinen nicht. Mit ${\cal O}(\overline R,G)$ bezeichnen wir die topologische multiplikative Gruppe aller auf $\overline R$ holomorphen Abbildungen mit Werten in $G$ ausgestattet mit der von der Supremumsnorm $\|.\|_{\overline R}$ erzeugten Topologie.

\begin{lmm}\label{hector-anita} Sei $R$ ein beliebiges Rechteck in $\mathbb C$. Dann lä\3t sich jede Abbildung $f\in{\cal O}(\overline R,G_1)$ in der Form $f=(1-h_1)\ldots(1-h_m)$, $m\geq 1$ schreiben mit $(1-h_j)\in{\cal O}(\overline R,G_1)$ und $\|h_j\|_{\overline R}<1$, $j=1,\ldots,m$.
\end{lmm}

\begin{proof}[{\bf Beweis}] 

\noindent Wir weisen nach, dass die Gruppe ${\cal O}(\overline R,G_1)$ zusammenhängend ist: Sei dazu $g$ eine beliebige Abbildungen aus ${\cal O}(\overline R,G_1)$ und $z_0$ der Mittelpunkt von $\overline R$. Es ist dann leicht zu sehen, dass die Abbildung $H:[0,1]\rightarrow{\cal O}(\overline R,G_1)$, welche durch die Formel $$(H(t))(z):=g(t(z-z_0)+z_0)\:,\ z\in\overline R$$ bestimmt ist und einen stetigen Weg in ${\cal O}(\overline R,G_1)$ darstellt, welcher das Element $H(1)=g$ mit dem Element $H(0)\equiv g(z_0)\in G_1$, also $H(0)\in{\cal O}(\overline R,G_1)$, verbindet. Da $g$ beliebig gewählt war und $G_1$ zusammenhängend ist, folgt daraus der Zusammenhang der topologischen Gruppe ${\cal O}(\overline R,G_1)$.

\noindent Wie schon im stetigen Fall, kann nun eine zusammenhängende topologische Gruppe durch eine offene Umgebung des Einselementes erzeugt werden. Speziell betrachten wir die Umgebung $Z:=\{(1-h)\in{\cal O}(\overline R,G_1)\mid \|h\|_{\overline R}<1\}$. Jedes Element $f\in{\cal O}(\overline R,G_1)$ lä\3t sich dann in der Form $$f=(1-h_1)\ldots(1-h_m)\:,\ m\geq 1\:$$ für endlich viele $(1-h_j)\in Z,\ j=1,\ldots,m$ schreiben, wobei der Nachweis analog zum ersten Punkt aus dem Beweis von Lemma \ref{horst-helena} ist.\footnote{Das in den Faktoren gewählte Plus- oder Minuszeichen ist dabei belanglos und kann je nach Kontext beliebig gewählt werden.}
\end{proof}

\begin{lmm}\label{hector-xandra} Es sei $\overline R\subset{\mathbb C}$ der Abschlu\3 eines Rechtecks. Dann kann man jede holomorphe Abbildung $f:\overline R\rightarrow A$ mit einem beliebigen Genauigkeitsgrad gleichmä\3ig bezüglich der Supremumsnorm $\|.\|_{\overline R}$ auf $\overline R$ durch Polynome der Art $$p(z)=\sum\limits^{N}_{j=0}z^jp_j\:,\ p_j\in A$$ approximieren.
\end{lmm}

\begin{bmr}\label{hector-victoria} Der Beweis ist angelehnt an den Beweis des verallgemeinerten Theorems von Runge (siehe \cite{gun}, Seite 36). Natürlich ist es für den Beweis nicht wichtig, abgeschlossene Rechtecke zu betrachten. Genausogut können zum Beispiel einfach zusammenhängende kompakte Mengen betrachtet werden.
\end{bmr}

\begin{proof}[{\bf Beweis}] Wir wählen eine hinreichend kleine offene Umgebung $U$ von $\overline R$, deren Rand aus einer rektifizierbaren Jordan-Kurve besteht, so dass die Abbildung $f$ noch holomorph auf $\overline U$ fortgesetzt werden kann. Die Cauchysche Integralformel liefert uns $$f(z)=\frac{1}{2\pi i}\int\limits_{\partial U}\frac{f(\zeta)}{\zeta-z}d\zeta\:,\ z\in \overline R\:.$$ Dieses Integral lä\3t sich auf $\overline R$ gleichmä\3ig durch Riemann-Summen mit beliebiger Genauigkeit approximieren: $$f(z)\approx\frac{1}{2\pi i}\sum\limits^{M}_{j=1}\frac{f(\zeta_j)}{\zeta_j-z}\Delta_j\:,\ \Delta_j:=\zeta_j-\zeta_{j-1}$$ für ein genügend gro\3es $M\in{\mathbb N}$. Es bleibt zu zeigen, dass für $\zeta\notin\overline R$ die komplexwertige Funktion $(\zeta-z)^{-1}$ auf $\overline R$ gleichmä\3ig durch Polynome approximiert werden kann. Sei dazu $\Gamma$ eine Kurve in ${\mathbb C}\setminus\overline R$, parametrisiert durch $\zeta:[0,\infty)\rightarrow{\mathbb C}$, so dass $|\zeta(t)|\stackrel{t\rightarrow\infty}{\longrightarrow}\infty$ gilt. Die Menge $$T:=\{t\in[0,\infty)\mid(\zeta(t)-z)^{-1}\ \text{ist beliebig genau durch Polynome approximierbar}\}$$ ist offenbar abgeschlossen, da eine gleichmä\3ig konvergente Folge durch Polynome approximierbarer Abbildungen selbst durch Polynome approximierbar ist. 

\newpage

\noindent $T$ ist nichtleer, da für $|\zeta(t_1)|>\sup\{|z|\mid z\in\overline R\}$ die Funktion $(\zeta(t_1)-z)^{-1}$ in einer offenen Kreisscheibe, welche $\overline R$ beinhaltet, holomorph ist und durch eine Potenzreihe repräsentiert werden kann. 

\begin{figure}[hbt] \begin{center}
\psfrag{Ga}{$\Gamma$}
\psfrag{zetat1}{$\zeta(t_1)$}
\psfrag{zeta}{$\zeta$}
\psfrag{z0}{}
\psfrag{z}{$z$}
\psfrag{Rq}{$\overline{R}$}
\psfrag{U}{$U$}
\includegraphics*[width=8cm]{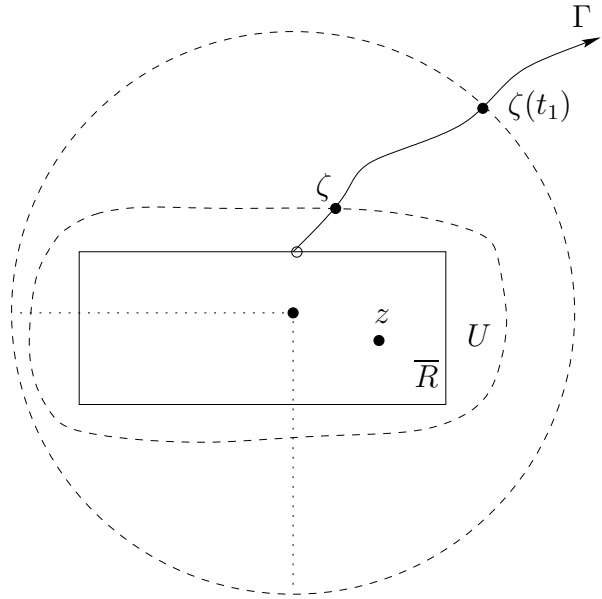}  
\end{center} \caption{\label{fig:ApproxF} Approximation der Funktion $\frac{1}{\zeta-z}$} \end{figure}

\noindent Schlie\3lich ist $T$ auch offen: Wir wählen $t_0\in T$ und zeigen, dass es eine Umgebung von $t_0$ gibt, die noch ganz in $T$ enthalten ist. Wir setzen 
\begin{equation}\label{hector-xandra-kroki}\varepsilon:=\inf\{|\zeta(t_0)-z|\mid z\in\overline R\}>0\:,
\end{equation}
und aufgrund der Stetigkeit von $\zeta$ als Funktion vom Parameter $t$ existiert ein $\delta>0$ derart, dass aus $|t_0-t|<\delta$ die Beziehung $|\zeta(t_0)-\zeta(t)|<\varepsilon$ folgt. Es bleibt zu zeigen, dass $t\in T$ für alle $t$ mit $|t_0-t|<\delta$ gilt. Wir können $(\zeta(t)-z)^{-1}$ als eine Potenzreihe in $(\zeta(t_0)-z)^{-1}$ schreiben, welche gleichmä\3ig auf $\overline R$ konvergiert:
 
\begin{eqnarray*}
\frac{1}{\zeta(t)-z} & = & \frac{1}{\zeta(t_0)-z+\zeta(t)-\zeta(t_0)}\\
 & = & \frac{1}{(\zeta(t_0)-z)\Big(1-\underbrace{\frac{\zeta(t_0)-\zeta(t)}{\zeta(t_0)-z}}_{\text{betragsmä\3ig}\ <1\ \text{wegen (\ref{hector-xandra-kroki})}}\Big)}\\
 & = & \frac{1}{\zeta(t_0)-z}\sum\limits^{\infty}_{k=0}\Big(\frac{\zeta(t_0)-\zeta(t)}{\zeta(t_0)-z}\Big)^{k}\ \text{mit Hilfe der geometrischen Reihe}\\
 & = & \sum\limits^{\infty}_{k=0}\frac{\big(\zeta(t_0)-\zeta(t)\big)^k}{\big(\zeta(t_0)-z\big)^{k+1}}\\
 & \approx & \sum\limits^{L}_{k=0}\frac{\big(\zeta(t_0)-\zeta(t)\big)^k}{\big(\zeta(t_0)-z\big)^{k+1}}\ \text{für ein genügend gro\3es}\ L\in{\mathbb N}\:.
\end{eqnarray*}
Da $t_0\in T$, ist auch $t\in T$, und da $T$ bezüglich der zusammenhängenden Menge $[0,\infty)$ offen, abgeschlossen und nichtleer ist, gilt $T=[0,\infty)$, und das Lemma ist bewiesen.
\end{proof}

\begin{lmm}[{\bf Approximationssatz}]\label{hector-nora} Es sei $f:\overline R\rightarrow G$ eine holomorphe Abbildung, wobei $\overline R\subset{\mathbb C}$ der Abschlu\3 eines Rechtecks ist. Dann existiert für jedes $\varepsilon>0$ eine holomorphe Abbildung $\tilde{f}:{\mathbb C}\rightarrow G$, so dass gilt: $$\|f-\tilde{f}\|_{\overline R}<\varepsilon\:,\ \text{beziehungsweise}\ \:\|1-f\tilde{f}^{-1}\|_{\overline R}<\varepsilon\:.$$
\end{lmm}

\begin{proof}[{\bf Beweis}] Ohne Beschränkung der Allgemeinheit kann man wie im Beweis von Satz \ref{horst-betty} annehmen, dass alle Bilder $f(z)$ zu der Zusammenhangskomponente $G_1$ gehören. Andernfalls kann man die Abbildung $(f(z_0))^{-1}f(z),\ z\in{\overline R}$ betrachten, wobei $z_0\in\overline R$ beliebig aber fest gewählt sei. Aus Lemma \ref{hector-anita} folgt $$f=(1-h_1)\ldots(1-h_m)=\exp({\ln(1-h_1)})\ldots\exp({\ln(1-h_m)})\:,\ m\geq 1\:,$$ wobei $(1-h_j)\in{\cal O}(\overline R,G_1)$ und $\|h_j\|_{\overline R}<1,\ j=1,\ldots,m$ ist. Für alle $j\in\{1,\ldots,m\}$ gilt: $$\ln(1-h_j)=\sum^{\infty}_{k=1}\frac{1}{k}h^k_j\:.$$ Aufgrund von Lemma \ref{hector-xandra} existieren Polynome $p_j:{\mathbb C}\rightarrow{A},\ j=1,\ldots,m$, welche die entsprechenden Abbildungen $\ln(1-h_j)$ auf $\overline R$ so genau gleichmä\3ig approximieren, dass die Abbildung $$\tilde{g}:=\exp p_1\ldots\exp p_m$$ das Gewünschte liefert. Damit ist das Lemma bewiesen.
\end{proof}

\noindent Nun haben wir alle Voraussetzungen, um folgenden Satz zu beweisen.

\begin{stz}\label{hector-manu} Seien $R_1$ und $R_2$ Rechtecke in $\mathbb{C}$, so dass ihre Vereinigung wieder ein Rechteck ist. Dann gibt es zu jeder Abbildung $f\in{\cal O}(\overline{R_1\cap R_2},G)$ Abbildungen $f_1\in\overline{\cal O}(R_1,G)$ und $f_2\in\overline{\cal O}(R_2,G)$ mit $f=f_1^{-1}f_2$ auf $R_1\cap R_2$.
\end{stz}

\begin{proof}[{\bf Beweis}] Sei $f\in{\cal O}(\overline{R_1\cap R_2},G)$ beliebig gewählt. Dann liefert uns der Approximationssatz mit $\varepsilon=\varepsilon_0$ (aus Lemma \ref{hector-amata}) eine holomorphe Abbildung $\tilde f:{\mathbb C}\rightarrow G$, für die $$f\tilde f^{-1}=1+g$$ auf $\overline{R_1\cap R_2}$ gilt, wobei $(1+g)\in\overline{\cal O}(R_1\cap R_2,G)$ und $\|g\|_{\overline{R_1\cap R_2}}<\varepsilon_0$ ist. Die Abbildung $(1+g)$ erfüllt die Voraussetzungen von Lemma \ref{hector-amata} und lä\3t sich in der Form $$1+g=(1+g_1)(1+g_2)$$ mit $(1+g_1)\in\overline{\cal O}(R_1,G)$ und $(1+g_2)\in\overline{\cal O}(R_2,G)$ schreiben (siehe die Setzung im Beweis von Lemma \ref{hector-amata}). Durch Umstellen erhalten wir $$f=\underbrace{(1+g_1)}_{=:f_1^{-1}}\underbrace{(1+g_2)\tilde f}_{=:f_2}$$ und der Satz ist bewiesen.
\end{proof}

\noindent Die Aussage des Satzes behält ihre Gültigkeit, wenn man Holomorphie nur im Inneren verlangt (siehe Satz \ref{hector-kata}).

\noindent Für eine beliebige nichtleere Menge $Y\subset{\mathbb C}$ sei ${\cal O}(Y,A)$ die Menge aller auf $Y$ holomorphen Abbildungen mit Werten in $A$, die offenbar einen Vektorraum über $\mathbb C$ bildet. Zur Definition einer Topologie betrachten wir ${\cal O}(Y,A)$ als Teilraum des Fréchetraumes $C^0(Y,A)$.

\begin{stz}[{\bf Faktorisierungssatz}]\label{hector-kata} Seien $R_1$ und $R_2$ Rechtecke in $\mathbb{C}$, so dass ihre Vereinigung wieder ein Rechteck ist. Dann gibt es zu jeder holomorphen Abbildung $f:R_1\cap R_2\rightarrow G$ holomorphe Abbildungen $f_1:R_1\rightarrow G$ und $f_2:R_2\rightarrow G$ mit $f=f_1^{-1}f_2$ auf $R_1\cap R_2$.
\end{stz}

\begin{proof}[{\bf Beweis}] Es seien $(R_{jn})_{n\in{\mathbb N}}$, $j=1,2$ Ausschöpfungen von $R_j$ durch Rechtecke, für die $R^n:=R_{1n}\cap R_{2n}$ und $R_{1n}\cup R_{2n}$ wieder Rechtecke sind, wobei $\overline R$$^n\subset R^{n+1}$ und $\bigcup\limits_{n\in {\mathbb N}}R^n=R_1\cap R_2$ gilt, und $R^1$ nichtleer ist (siehe Abbildung \ref{fig:bild9}). 

\begin{figure}[htb]
\begin{center}
\psfrag{R1n}{\hspace{-2mm}$R_{1n}$}
\psfrag{R2n}{$R_{2n}$}
\psfrag{Rn}{$R^n$}
\psfrag{R1}{$R_1$}
\psfrag{R2}{$R_2$}
  \includegraphics*[width=7cm]{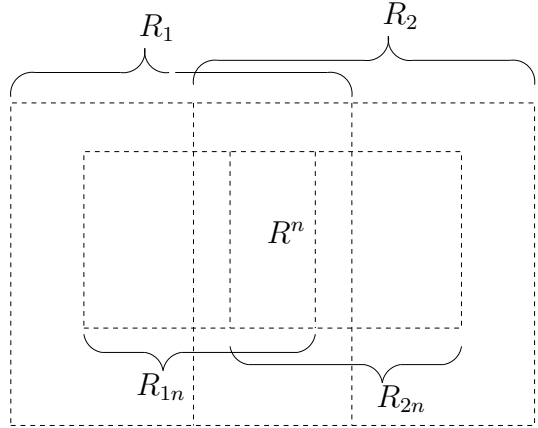}
 \end{center} 
 \caption{\label{fig:bild9} Ausschöpfungen der Rechtecke $R_1$ und $R_2$}
\end{figure}

\noindent Gegeben sei die in $R_1\cap R_2$ holomorphe Abbildung $f$. Nach Satz \ref{hector-manu} gibt es für jedes $n\in {\mathbb N}$ auf $\overline R_{jn}$ stetige und in $R_{jn}$ holomorphe Abbildungen $f_{jn}$, $j=1,2$ mit $f=f_{1n}^{-1}f_{2n}$ auf $R^n$. Es sei $$v_n:=f_{jn}f_{j,n+1}^{-1}$$ auf $R_{jn}$, wobei auf $R^n$ gilt $f_{1n}^{-1}f_{2n}=f_{1,n+1}^{-1}f_{2,n+1}$, also $f_{1n}f_{1,n+1}^{-1}=f_{2n}f_{2,n+1}^{-1}$, und damit ist $v_n$ wohldefiniert, holomorph in $R_{1n}\cup R_{2n}$ und stetig auf $\overline{R_{1n}\cup R_{2n}}$. Durch wiederholte Anwendung des Approximationssatzes erhält man holomorphe Abbildungen $g_n:{\mathbb C}\rightarrow G$ mit $$\|1-(g_nv_n)g_{n+1}^{-1}\|_{\overline{R_{1n}\cup R_{2n}}}<\frac{1}{2^{n+1}}\:,\ n\in {\mathbb N}\:,$$ wobei $g_1:=1$ gesetzt wird. Dann ist $$h_n:=\prod\limits^{\infty}_{k=n}(g_{k}v_{k}g_{k+1}^{-1})$$ auf $\overline{R_{1n}\cup R_{2n}}$ gleichmä\3ig konvergent, so dass $h_n$ in $R_{1n}\cup R_{2n}$ holomorph ist, und stetig auf $\overline{R_{1n}\cup R_{2n}}$. Es gilt: $$\|1-h_n\|_{\overline{R_{1n}\cup R_{2n}}}<\frac{1}{2^n}\e^{\frac{1}{2^n}}<1\:,$$ also $h_n\in{\cal O}(R_{1n}\cup R_{2n},G)$. Wegen $h_n=(g_nv_ng_{n+1}^{-1})h_{n+1}$ auf $R_{1n}\cup R_{2n}$ gilt $g_n^{-1}h_n=v_ng_{n+1}^{-1}h_{n+1}$ und es folgt $f_{jn}^{-1}g_n^{-1}h_n=f_{j,n+1}^{-1}g_{n+1}^{-1}h_{n+1}$ auf $R_{jn}$ für $j=1,2$. Setzen wir $$f_j:=f_{jn}^{-1}g_n^{-1}h_n\:,\ j=1,2$$ auf $R_{jn}$, so ist $f_j:R_j\rightarrow G$ wohldefiniert und holomorph. Dann gilt: $$f_1^{-1}f_2=f_{1n}^{-1}g_n^{-1}h_nh_n^{-1}g_nf_{2n}=f_{1n}^{-1}f_{2n}=f$$ auf $R^n$ für alle $n\in{\mathbb N}$, und folglich auch auf $R_1\cap R_2$.
\end{proof}

\noindent Wie im stetigen Fall erhält man den folgenden Satz.

\begin{stz}\label{hector-michaela} Es sei ${\cal O}^G$ die Garbe der Keime holomorpher Abbildungen mit Werten in $G$, und $R_1,\ldots,R_n$, $n\in{\mathbb N}$ seien Rechtecke in $\mathbb C$, so dass die Durchschnitte $(R_1\cup\ldots\cup R_k)\cap R_{k+1}$ und die Vereinigungen $R_1\cup\ldots\cup R_{k+1}$, $k=1,\ldots,n-1$ wieder Rechtecke sind. 

\noindent Dann gilt $H^1({\cal R},{\cal O}^G)=1$ für die Überdeckung ${\cal R}:=(R_k)_{1\leq k\leq n}$ von $R_1\cup\ldots\cup R_n$.
\end{stz}

\begin{proof}[{\bf Beweis}] Der Nachweis erfolgt analog zum Beweis von Satz \ref{horst-petra} im stetigen Fall:

\noindent Es sei ${\cal F}={\cal O}^G$. Für $n=1$ ist die Aussage trivial. Für $n>1$ sind wegen Satz \ref{hector-kata} $(R_1\cup\ldots\cup R_k,R_{k+1})$ für alle $k\in\{1,\ldots,n-1\}$ $\cal F$-Paare und damit ist $(R_1,\ldots,R_n)$ eine $\cal F$-Kette, woraus mit Satz \ref{harry-kathi} folgt, dass $H^1({\cal R},{\cal O}^G)$ trivial ist.
\end{proof}

\begin{stz}\label{hector-sandra} Sei $R$ ein Rechteck in $\mathbb C$ und ${\cal R}=(R_i^j)_{1\leq i\leq n_j,1\leq j\leq m}$ mit $n_j,m\in{\mathbb N}$ eine Überdeckung von $R$, wobei die $R_i^j$ Rechtecke in $\mathbb C$ sind, für die gilt:
\begin{enumerate}
	\item für alle $j\in\{1,\ldots,m\}$ sind die Durchschnitte $(R_1^j\cup\ldots\cup R_k^j)\cap R_{k+1}^j$ und die Vereinigungen $R_1^j\cup\ldots\cup R_{k+1}^j$, $k=1,\ldots,n_j-1$ Rechtecke in $\mathbb C$ und
	\item mit $R^j:=\bigcup\limits_{1\leq i\leq n_j}R_i^j$ für alle $j\in\{1,\ldots,m\}$ sind $(R^1\cup\ldots\cup R^k)\cap R^{k+1}$ und $R^1\cup\ldots\cup R^{k+1}$, $k=1,\ldots,m-1$ wieder Rechtecke in $\mathbb C$.
\end{enumerate}
Dann ist $H^1({\cal R},{\cal O}^G)=1$.
\end{stz}

\begin{figure}[htb]
\begin{center}
   \includegraphics*[width=14cm]{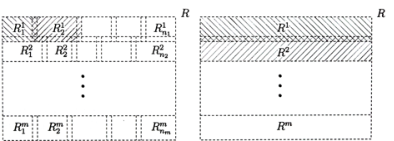}
 \end{center} 
 \caption{\label{fig:bild8} Überdeckung von $R$}
\end{figure}

\begin{proof}[{\bf Beweis}] Hier ist wieder ${\cal F}={\cal O}^G$. Wie im Beweis von Satz \ref{hector-michaela} erhalten wir sofort für alle $j\in\{1,\ldots,m\}$ $\cal F$-Ketten ${\cal R}^j:=(R_1^j,\ldots,R_{n_j}^j)$ und eine weitere $\cal F$-Kette $(R^1,\ldots,R^m)$. Laut Definition handelt es sich bei ${\cal R}={\cal R}^1\cup\ldots\cup{\cal R}^m$ um ein $\cal F$-Feld und Satz \ref{harry-conny} liefert uns wie gewünscht $H^1({\cal R},{\cal O}^G)=1$.
\end{proof}

\noindent Es geht nun darum, diesen Satz zu verallgemeinern, d.h. wir wollen beliebige offene Überdeckungen betrachten, und dabei nicht nur Rechtecke überdecken. 

\begin{stz}[{\bf Zerfällungssatz}]\label{hector-gabi} Sei $Y$ ein einfach zusammenhängendes Gebiet in $\mathbb C$, das mehr als einen Randpunkt hat. Dann ist $H^1(Y,{\cal O}^G)=1$.
\end{stz}

\begin{proof}[{\bf Beweis}] Mit Hilfe des Riemannschen Abbildungssatzes bildet man $Y$ eineindeutig und konform auf ein Rechteck $R$ in $\mathbb C$ ab, indem man beide Gebiete zunächst auf das Innere des Einheitskreises abbildet. Sei $\varphi:Y\rightarrow R$ eine solche Abbildung. Dabei vermittelt die bijektive holomorphe Abbildungsfunktion $\varphi$ eine bijektive Abbildung der holomorphen Garbe ${\cal O}^G$ auf sich, d.h. wir betrachten die Keime der Abbildungen nicht bezüglich der Punkte $a\in Y$, sondern bezüglich ihrer Bildpunkte $\varphi(a)\in R$, und es genügt zu zeigen $H^1(R,{\cal O}^G)=1$.

\noindent Sei ${\cal U}=(U_k)_{k\in K}$ eine beliebige Überdeckung mit offenen Mengen von $Y$. Wir betrachten für $n=1,2,\ldots$ die Abschlüsse der Rechtecke $$R_n:=\Big\{z\in R\mid\dist(z,\partial R)>\frac{1}{C\cdot n}\Big\}\:,$$ wobei die Konstante $C$ genügend gro\3 gewählt sei, damit $R_1$ nichtleer ist. Dann gilt $\bigcup\limits^{\infty}_{n=1}\overline R_n=\bigcup\limits^{\infty}_{n=1} R_n=R$. Da $\overline R_n$ kompakt ist, gibt es endlich viele Mengen $U_{k_l},\ {1\leq l\leq m_n}$ der Überdeckung ${\cal U}=(U_k)_{k\in K}$, für die $\overline R_n\subset\bigcup\limits_{l=1}^{m_n}U_{k_l}$ gilt. Dann ist ${\cal U}_n:=(U_{k_l}\cap R_n)_{1\leq l\leq m_n}$ eine endliche Überdeckung von $R_n$, so dass eine Verfeinerung ${\cal R}_n=(R_i^j)_{1\leq i\leq n_j,1\leq j\leq m}$ existiert, welche den Bedingungen von Satz \ref{hector-sandra} genügt. Damit gilt $H^1({\cal R}_n,{\cal O}^G)=1$, und weiterhin $H^1({\cal U}|_{R_n},{\cal O}^G)=1$ wegen Korollar \ref{harry-clothilde}, also $$H^1(R_n,{\cal O}^G)=1\ \forall n\in{\mathbb N}$$ da $\cal U$ beliebig gewählt war. Wir betrachten an dieser Stelle die Überdeckung $${\cal R}:=(R_n)_{n\in{\mathbb N}}$$ von $R$. Dann gilt für alle $\Psi\in H^1(R,{\cal O}^G)$ nun $\Psi|_{R_n}=1\ \forall n\in{\mathbb N},$ und wir können Satz \ref{harry-lucia} anwenden. Dieser liefert uns ein $\Phi\in H^1({\cal R},{\cal O}^G)$ mit $[\Phi]=\Psi$. Es sei $f\in{\cal Z}^1({\cal R},{\cal O}^G)$ ein beliebiger $1$-Kozyklus mit $[f]=\Phi$, und wir zeigen, dass ein $w\in {\cal C}^0({\cal R},{\cal O}^G)$ mit $f_{n,n+1}=w_n^{-1}w_{n+1}$ auf $R_n$ für alle $n\in{\mathbb N}$ existiert:

\noindent Durch wiederholte Anwendung des Approximationssatzes erhalten wir holomorphe Abbildungen $g_n:{\mathbb C}\rightarrow G$ mit $$\|1-(g_nf_{n,n+1})g_{n+1}^{-1}\|_{\overline R_{n-1}}<\frac{1}{2^{n+1}},\ n\in{\mathbb N}\:,$$ wobei $g_1:=1$ gesetzt wird. Dann ist $$h_n:=\prod^{\infty}_{k=n}(g_kf_{k,k+1}g_{k+1}^{-1})$$ auf $\overline R_{n-1}$ gleichmä\3ig konvergent, so dass $h_n$ auf $R_{n-1}$ holomorph ist. Weiterhin gilt: $$\|1-h_n\|_{\overline R_{n-1}}<\frac{1}{2^n}\e^{\frac{1}{2^n}}<1\:,$$ so dass $h_n(z)$ für alle $z\in\overline R_{n-1}$invertierbar ist, also $h_n|_{R_{n-1}}\in{\cal O}^G(R_{n-1})$ gilt. Wegen $$h_n=(g_nf_{n,n+1}g_{n+1}^{-1})h_{n+1}$$ auf $R_n$ folgt $h_n\in{\cal O}^G(R_{n})$, denn der erste Faktor ist Produkt dreier dort definierter holomorpher Abbildungen mit invertierbaren Werten und der zweite Faktor $h_{n+1}$ ist nach dem eben Gezeigten ebenfalls holomorph mit invertierbaren Werten auf $R_n$. Dann gilt $f_{n,n+1}=(h_n^{-1}g_n)^{-1}h_{n+1}^{-1}g_{n+1}$ auf $R_n$. Setzen wir $$w_n:=h_n^{-1}g_n\ \:\text{und}\ \:w_{n+1}:=h_{n+1}^{-1}g_{n+1}$$ auf $R_n$ beziehungsweise auf $R_{n+1}$, so sind $w_n:R_n\rightarrow G$ und $w_{n+1}:R_{n+1}\rightarrow G$ holomorph und es gilt auf $R_n\cap R_{n+1}=R_n$: $$f_{n,n+1}=w_n^{-1}w_{n+1}\:.$$
\end{proof}

\newpage

\printindex

\newpage

\section*{Selbstständigkeitserklärung}

\noindent Ich erkläre, dass ich die vorliegende Arbeit selbstständig und nur unter Verwendung der angegebenen Literatur und Hilfsmittel angefertigt habe.

\noindent $ $ 

\noindent Berlin, den 9. September 2004

\noindent $ $ 

\noindent Katrin Kaden

\newpage

\section*{Thesen}

\begin{itemize}
	\item Seien $\cal G$ eine Garbe über einem topologischen Raum $X$, $\cal U$ eine offene Überdeckung von $X$ und $\Psi$ ein Element aus der Kohomologiemenge $H^1(X,{\cal G})$, dessen Einschränkungen auf die Mengen von $\cal U$ jeweils die neutralen Elemente sind. Dann existiert ein $\Phi$ aus der Kohomologiemenge $H^1({\cal U},{\cal G})$, dessen Äquivalenzklasse mit $\Psi$ übereinstimmt (Satz \ref{harry-lucia}).
	\item Seien $\cal F$ eine Garbe über $X$ und $\cal U$ eine $\cal F$-Kette beziehungsweise ein $\cal F$-Feld. Dann ist die Kohomologiemenge $H^1({\cal U},{\cal F})$ trivial (Sätze \ref{harry-kathi} und \ref{harry-conny}).
  \item Seien $X$ ein metrischer Raum, der sich als Vereinigung abzählbar vieler kompakter Teilmengen darstellen lä\3t, $G$ die Menge der invertierbaren Elemente einer Banachalgebra und $G_1\subset G$ diejenige Zusammenhangskomponente, welche das Einselement enthält. Dann gelten:	
\end{itemize}

\begin{tabular}{p{2,7cm}|p{5,1cm}|p{5,1cm}}
 & stetiger Fall & holomorpher Fall \\ \hline Fortsetzungs- bzw. Approximationssatz & Seien $W\subset X$ abgeschlossen und kontrahierbar sowie $f:W\rightarrow G_1$ eine stetige Abbildung. Dann existiert eine stetige Abbildung $\tilde f:X\rightarrow G$ aus der Zusammenhangskomponente der Einsabbildung, die auf $W$ mit $f$ übereinstimmt (Lemma \ref{horst-karla}). & Seien $\overline R\subset{\mathbb C}$ der Abschlu\3 eines Rechtecks und $f:\overline R\rightarrow G$ eine holomorphe Abbildung. Dann gibt es eine holomorphe Abbildung $\tilde f:{\mathbb C}\rightarrow G$, welche die Abbildung $f$ auf $\overline R$ beliebig genau approximiert (Lemma \ref{hector-nora}). \\ \hline Faktorisierungs-sätze & Seien $U_1,U_2$ offene Teilmengen von $X$ mit einem nichtleeren kontrahierbaren Durchschnitt. Dann lä\3t sich jede stetige Abbildung $f:U_1\cap U_2\rightarrow G$ faktorisieren (Satz \ref{horst-betty}). & Seien $R_1,R_2$ Rechtecke in ${\mathbb C}$, deren Vereinigung wieder ein Rechteck ist. Dann lä\3t sich jede holomorphe Abbildung $f:R_1\cap R_2\rightarrow G$ faktorisieren (Satz \ref{hector-kata}). \\ \hline Zerfällungs-sätze & Seien ${\cal C}^G$ die Garbe der Keime stetiger Abbildungen über $X$ mit Werten in $G$ und $U_1,\ldots,U_n,\ n\in{\mathbb N}$ offene Mengen in $X$, wobei die Durchschnitte von $U_{k+1}$ und $U_1\cup\ldots\cup U_k,\ k=1,\ldots,n-1$ kontrahierbar sind. Dann ist die Kohomologiemenge $H^1((U_k)_{1\leq k\leq n},{\cal C}^G)$ trivial (Satz \ref{horst-petra}). & Seien ${\cal O}^G$ die Garbe der Keime holomorpher Abbildungen über ${\mathbb C}$ mit Werten in $G$ und $Y\subset{\mathbb C}$ ein einfach zusammenhängendes Gebiet, welches mehr als einen Randpunkt hat. Dann ist die Kohomologiemenge $H^1(Y,{\cal O}^G)$ trivial (Satz \ref{hector-gabi}).
\end{tabular}

\end{document}